\RequirePackage[l2tabu, orthodox]{nag} 
\documentclass{amsart}
\usepackage[abbrev]{amsrefs}
\usepackage{yhmath}
\usepackage{todonotes}
\usepackage{bm}
\usepackage{dsfont}

\usepackage{mathtools}
\mathtoolsset{showonlyrefs=true}
\usepackage{amssymb}
\usepackage{amsthm}
\usepackage{mathrsfs}
\usepackage{enumerate}
\usepackage[shortlabels]{enumitem}
\usepackage[title]{appendix}

\usepackage[T1]{fontenc} 
\usepackage[top=1in, bottom=1in, left=1in, right=1in]{geometry}
\usepackage[numbered,autolinebreaks,useliterate,bw]{mcode}
\newcommand{\argmin}{\operatorname*{argmin}}

\def\spt{\mathop{\mathrm{spt}}\nolimits}

\theoremstyle{plain}
 \newtheorem{theorem}{Theorem}[section]
 \newtheorem{lemma}[theorem]{Lemma}
 
 \newtheorem{proposition}[theorem]{Proposition}

\theoremstyle{definition}
 \newtheorem{definition}[theorem]{Definition}
 
\newtheorem{assumption}[theorem]{Assumption}
 
\theoremstyle{remark}
 \newtheorem{remark}[theorem]{Remark}

\numberwithin{equation}{section}

\begin{document}
 \title[Error estimate for regularized optimal transport problems via Bregman divergence]{Error estimate for regularized optimal transport problems \\via Bregman divergence}

 \author[K.~Morikuni]{Keiichi Morikuni}
 \address[K.~Morikuni]{Institute of Systems and Information Engineering, University of Tsukuba,
 1-1-1 Tennodai, Tsukuba-shi, Ibaraki 305-8573, Japan}
 \email{morikuni@cs.tsukuba.ac.jp}
 \author[K.~Sakakibara]{Koya Sakakibara}
 \address[K.~Sakakibara]{Faculty of Mathematics and Physics, Institute of Science and Engineering, Kanazawa University,
 Kakuma-machi, Kanazawa-shi, Ishikawa 920-1192, Japan;
 RIKEN iTHEMS,
 2-1 Hirosawa, Wako-shi, Saitama 351-0198, Japan}
 \email{ksakaki@se.kanazawa-u.ac.jp}
 \author[A.~Takatsu]{Asuka Takatsu}
 \address[A.~Takatsu]{Graduate School of Mathematical Sciences, The University of Tokyo, 
 3-8-1 Komaba Meguro-ku Tokyo 153-8914, Japan;
 RIKEN Center for Advanced Intelligence Project (AIP), 
Nihonbashi 1-chome Mitsui Building, 15th floor,
1-4-1 Nihonbashi, Chuo-ku, Tokyo 103-0027, Japan.}
 \email{asuka-takatsu@g.ecc.u-tokyo.ac.jp}
 \begin{abstract}
Optimal transport problems
on a finite set can be efficiently and approximately solved by using entropic regularization. 
This study proposes using Bregman divergence regularization for optimal transport problems. 
We discuss the required properties for Bregman divergence, provide a non-asymptotic error estimate for the regularized problem, and show that the error estimate becomes faster than exponential.
Numerical experiments show that our regularization outperforms that using the Kullback--Leibler divergence.
\end{abstract}

 \maketitle

\section{Introduction}\label{1}
\emph{Optimal transport theory} allows for measuring the difference between two probability measures.
Innumerable applications of optimal transport theory include mathematics, physics, economics, statistics, computer science, and machine learning.
This work focuses on  optimal transport theory on  finite sets.

For $K\in \mathbb{N}$, define 
\[
\mathcal{P}_K:=\left\{ z=(z_k)\in \mathbb{R}^K \biggm| z_k\geq 0 \ \text{for any } k,\ 
\sum_{k} z_k=1\right\}.
\]
Here and hereafter, $k$ runs over $1,2,\ldots, K$.
Fix $I,J \in \mathbb{N}$.
Unless indicated otherwise, 
$i$ and $j$ run over $1,2,\ldots, I$ and $1,2,\ldots, J$, respectively.
For $x\in \mathcal{P}_I$ and $y\in \mathcal{P}_J$,
define $x\otimes y \in \mathcal{P}_{I \times J}$ by 
\[
(x\otimes y)_{ij}:=x_iy_j,
\]
and set 
 \begin{align}
 \Pi(x,y):=
 \left\{
 \Pi=(\pi_{ij})
 \in\mathcal{P}_{I \times J}
 \Biggm| 
 \sum_{l=1}^{J}\pi_{il}=x_i,
 \sum_{l=1}^{I}\pi_{lj}=y_j 
 \ \text{for any } i,j
 \right\},
 \end{align}
where we identify $\mathcal{P}_{I\times J}$ with a subset of $\mathbb{R}^{I\times J}$.
An element in $\Pi(x,y)$ is called a \emph{transport plan} between $x$ and $y$.
Note that $\Pi(x,y)$ is a compact set, in particular, a convex polytope, and contains $x\otimes y$.
Let us fix $C=(c^{ij})\in\mathbb{R}^{I \times J}$ and 
define a map~$\langle C, \cdot \rangle :\mathcal{P}_{I \times J} \to \mathbb{R}$ by
\[
 \langle C,\Pi\rangle:=\sum_{i,j} c^{ij}\pi_{ij}.
\]
Consider a linear program of the form
 \begin{align}\label{eq:LP}
\inf_{\Pi\in\Pi(x,y)}\langle C,\Pi\rangle,
 \end{align}
commonly referred to as an \emph{optimal transport problem}.
Since the function~$\langle C, \cdot \rangle $ is linear, in particular continuous on a compact set~$\Pi(x,y)$,
the problem~\eqref{eq:LP} always admits a minimizer, which is not necessarily unique.
A minimizer of the problem~\eqref{eq:LP} is called an \emph{optimal transport plan} between $x$ and~$y$.

This paper considers a regularized optimal transport problem via Bregman divergence,
which is a generalization of the Kullback--Leibler divergence using a strictly convex function.
\begin{definition}
Let $U$ be a continuous, strictly convex function on $[0,1]$ with $U\in C^1((0,1])$.
For $z, w\in \mathcal{P}_K$, the \emph{Bregman divergence}
associated with $U$ of $z$ with respect to $w$ is given by
\[
D_U(z,w):=\sum_{k} d_U(z_k, w_k),
\]
where $d_U:[0,1]\times (0,1] \to \mathbb{R}$ is defined for $r\in [0,1]$ and $r_0\in (0,1]$ by 
\[
d_U(r,r_0):=U(r)-U(r_0)-(r-r_0)U'(r_0)
\]
and is naturally extended as 
a function on $[0,1] \times [0,1]$
valued in $[0,\infty]$
(see Lemma~\ref{behaviorU}).
\end{definition}

For example, the Bregman divergence associated with $U(r)=r\log r$ coincides with the Kullback--Leibler divergence.

Let us consider a regularized problem of the form
 \begin{align}
\inf_{\Pi\in\Pi(x,y)} \left(\langle C,\Pi\rangle+\varepsilon D_U(\Pi,x \otimes y)\right)
 \quad
 \text{for }\varepsilon>0.
 \label{eq:OTe}
 \end{align}
By the continuity and strict convexity of $U$, $D_U(\cdot, x\otimes y)$ is continuous and strictly convex on a convex polytope~$\Pi(x,y)$.
Consequently, the problem~\eqref{eq:OTe} always admits a unique minimizer,
denoted by $\Pi^{U}(C,x,y,\varepsilon)$. 
Then,
\begin{align}\label{conv}
\lim_{\varepsilon \downarrow 0}
\langle C,\Pi^{U}(C,x,y,\varepsilon)\rangle
=
\inf_{\Pi\in\Pi(x,y)}\langle C,\Pi\rangle
\end{align}
holds (see Subsection~\ref{error}).
To obtain a quantitative error estimate of \eqref{conv},
we require the following two assumptions.
See Subsections~\ref{subxy} and \ref{verify} to verify that the assumptions are reasonable.
\begin{assumption}\label{xy}
$\Pi(x,y)\neq \argmin_{\Pi\in\Pi(x,y)}\langle C,\Pi\rangle$.
\end{assumption}

\begin{assumption}\label{assumption}
Let $U\in C([0,1]) \cap C^1((0,1])\cap C^2((0,1))$ satisfy 
$U''>0$ on~$(0,1)$ and $\lim_{h\downarrow0}U'(h)=-\infty$.
In addition, $r\mapsto r U''(r)$ is non-decreasing in~$(0,1)$.
\end{assumption}

We introduce notions to describe our quantitative error estimate of \eqref{conv}.
\begin{definition}\label{def:DU}
Let $U$ be a continuous, strictly convex function on $[0,1]$ with $U\in C^1((0,1])$.
Define $\mathfrak{D}_U(x,y)$ for $x\in\mathcal{P}_I$ and $y\in\mathcal{P}_J$ 
by
\begin{align}
 \mathfrak{D}_U(x,y)
\coloneqq\sup_{\Pi\in\Pi(x,y)}D_U(\Pi,x\otimes y).
 \label{eq:def_entropic_radius}
 \end{align}
\end{definition} 
\begin{definition}\label{def:sg}
The \emph{suboptimality gap}
of $x\in\mathcal{P}_I$  and $y\in\mathcal{P}_J$
with respect to $C\in \mathbb{R}^{I \times J}$ is defined by
\begin{align}
 \Delta_C(x,y)
 \coloneqq
 \inf_{V' \in V(x,y) \setminus \argmin_{V\in V(x,y)} \langle C,V\rangle }\langle C,V'\rangle
 -
 \inf_{V\in V(x,y)} \langle C, V \rangle,
 \end{align}
where $V(x,y)$ is the set of vertices of $\Pi(x,y)$ and set $\inf \emptyset:=\infty$. 
\end{definition}
In Subsection~\ref{subxy}, we verify $\mathfrak{D}_U(x,y),\Delta_C(x,y)\in (0,\infty)$ under Assumption~\ref{xy}.
We also confirm in Subsection~\ref{welldefined} that Definition~\ref{def:rU} below is well-defined.

\begin{definition}\label{def:rU}
Under Assumption~\ref{assumption},
let $e_U$ denote
the inverse function of $U':(0,1]\to U'((0,1])$.
For $x\in \mathcal{P}_I$ and $y\in \mathcal{P}_J$,
let $R_U(x,y)\in [1/2,1)$ satisfy
\[
U'(R_U(x,y))-U'(1-R_U(x,y))=\mathfrak{D}_U(x,y),
\]
which is uniquely determined.
Define $\nu_U(x,y)\in \mathbb{R}$ by 
\begin{align}\label{eq:delta_U}
\nu_U(x,y):=\sup_{r \in (0,R_U(x,y)]} \left(U'(1-r)+r U''(r)\right).
\end{align}
\end{definition}
Our main result is as follows.
\begin{theorem} \label{thm:convergence}
Under~Assumptions~{$\ref{xy}$ and $\ref{assumption}$}, 
the interval 
\[
\left(0,\frac{\Delta_C(x,y) R_U(x,y)}{\mathfrak{D}_U(x,y)}\right] \cap
\left(0, \frac{\Delta_C(x,y)}{\mathfrak{D}_U(x,y)+\nu_U(x,y)-U'(1)}\right]
\]
is well-defined and nonempty.
In addition,
\begin{align*}
&\langle C,\Pi^{U}(C,x,y,\varepsilon)\rangle
-\inf_{\Pi\in \Pi(x,y)}\langle C,\Pi\rangle\\
&\leq 
\Delta_C(x,y) \cdot 
e_U\left(-\frac{\Delta_C(x,y)}{\varepsilon}+\mathfrak{D}_U(x,y)
+\nu_U(x,y)\right)
 \end{align*}
holds for $\varepsilon$ in the above interval.
\end{theorem}

Before explaining the merit of our theorem,
let us review related results.
Computing an exact solution of a large-scale optimal transport problem becomes problematic when, say,  $N := \max\{I, J\} > 10^4$.
The best-known practical complexity~$\widetilde{O}(N^3)$ is attained by an interior point algorithm in~\cite{PW09}*{Section~5}, where $\widetilde{O}$ omits polylogarithmic factors.
Though Chen et al.~\cite{CKLPGS22FOCS}*{Informal Theorem~I.3} improve this complexity to $(N^2)^{1+o(1)}$ and Jambulapati et al.~\cite{JST19NeurIPS}*{Theorem~2.4} provide an algorithm that finds an $\epsilon$-approximation in $\widetilde{O}(N^2 / \epsilon)$, their practical implementations have not been developed.

The tractability of the problem~\eqref{eq:LP} is improved by introducing entropic regularization to its objective function, that is,
$\inf_{\Pi\in\Pi(x,y)}\left(\langle C,\Pi\rangle-\varepsilon S(\Pi)\right)$,
where 
\begin{equation}
    S(z):=-\sum_{k} z_k \log z_k, \quad z \in \mathcal{P}_K
    \label{eq:Shannon}
\end{equation}
is the \emph{Shannon entropy}.
Here, we put $0 \log 0:=0$ due to the continuity 
$ \lim_{r\downarrow0} r\log r=0$. 
Fang~\cite{Fan92} introduces the Shannon entropy to regularize generic linear programs.
By the continuity and the strict convexity of the Shannon entropy,
the entropic regularized problem always has a unique minimizer for each value of the regularization parameter.
Cominetti and San Mart\'{i}n~\cite{CM94}*{Theorem~5.8} prove that the minimizer of the regularized problem converges exponentially to a certain minimizer of the given problem as the regularization parameter goes to zero.
Weed~\cite{Wee18COLT}*{Theorem~5} provides a quantitative error estimate of the regularized problem, whose convergence rate is exponential.  
Note that
the entropic regularization allows us to develop approximation algorithms for the problem~\eqref{eq:LP}.
We refer to \cite{PeyreCuturi2019} and references therein.

Different types of regularizers have been introduced in recent studies.
For example, Muzellec et al.~\cite{MNPN17AAAI} use the Tsallis entropy for ecological inference.
Dessein et al.~\cite{DPR18} and Daniels et al.~\cite{DMH21} introduce the Bregman and $f$-divergences, respectively, to regularize optimal transport problems.
Apart from entropy and divergence, Klatt et al.~\cite{KTM20SIMODS} use convex functions of Legendre type for regularization.
These studies discuss convergence when the regularization parameter goes to zero.
In particular, Daniels et al.~\cite{DMH21}*{Theorem~4.3} provide a quantitative error estimate when the regularizer is an $f$-divergence.

Regularization by the Shannon entropy is the same as that by the Kullback--Leibler divergence.
Here, the \emph{Kullback--Leibler divergence} of $z \in \mathcal{P}_K$ with respect to $w \in \mathcal{P}_K$ is given by
\begin{align*}
 D_{\mathrm{KL}}(z,w)
 :=\sum_{k }^K z_k \left(\log z_k-\log w_k\right).
 \end{align*}
Note that the Kullback--Leibler divergence and its dual are the unique members that belong to both the Bregman and $f$-divergence classes (see \cite{Ama09IEEE} for instance).
Define $U_o\in C([0,\infty))\cap C^{\infty}((0,\infty))$ by 
\begin{align}
U_o(r):=
\begin{cases}
r \log r &\text{for }r\in (0,\infty),\\
0& \text{for }r=0.
\end{cases}
\label{eq:U_log}
\end{align}
Then, $D_{U_o}=D_{\mathrm{KL}}$ holds on $\mathcal{P}_K\times \mathcal{P}_K$.
There are other strictly convex functions~$U$ such that $D_U=D_{\mathrm{KL}}$ holds (see Subsection~\ref{UtoL}).

In this paper, we find a class of Bregman divergences 
by which the regularization performs better than that by the Kullback--Leibler divergence. 
Moreover, our main result Theorem~\ref{thm:convergence} with the case $U=U_o$
recovers Weed's result~\cite{Wee18COLT}*{Theorem~5}.
Theorem~\ref{thm:convergence} with 
the relation~\eqref{eu} guarantees that  
the regularized optimal value approaches the true optimal value faster than exponential
(see Subsection~\ref{verify}).
Numerical experiments demonstrate that a Bregman divergence gives smaller errors than the Kullback--Leibler divergence.

 This paper is organized as follows. 
 In Section \ref{sec:preliminaries}, 
 we verify that Assumptions~\ref{xy} and~\ref{assumption} are reasonable.
Section~\ref{sec:error_estimate} proves Theorem~\ref{thm:convergence}.
In Section~\ref{normalization}, 
we show that 
the normalization of $U$ does not affect the error estimate in Theorem~\ref{thm:convergence}.
We then consider the effect of scaling of data and the domain of $U$ on the error estimate.
Section~\ref{example} provides examples of $U$ satisfying Assumption~\ref{assumption}.
In Section \ref{sec:numerical_experiments},
numerical experiments show that regularization using a certain family of Bregman divergences outperforms that using the Kullback--Leibler divergence.
Finally, in Section~\ref{sec:condluding_remarks}, we summarize the contents of this paper and give directions for future research.

 \section{Preliminaries}
 \label{sec:preliminaries} 
In this section, 
we verify that Assumptions~\ref{xy}, \ref{assumption} are reasonable 
and Definition~\ref{def:rU} is well-defined.
We also show that $\mathfrak{D}_U(x,y), \Delta_{C}(x,y)\in (0,\infty)$ under 
Assumption~\ref{xy}. 
Throughout, as in the introduction, we fix $I, J \in\mathbb{N}$ and take
$C\in \mathbb{R}^{I\times J}$, 
$x\in \mathcal{P}_I$, and
$y\in \mathcal{P}_J$.
Let 
\[
\Omega:=\mathbb{R}^{I\times J}\times \mathcal{P}_I\times \mathcal{P}_J\times (0,\infty)
\]
and $U$ be a continuous, strictly convex function on~$[0,1]$ with $U\in C^1((0,1])$, unless otherwise stated.

By the strict convexity of $U$ on $[0,1]$,
\[
d_U(r,r_0):=U(r)-U(r_0)-(r-r_0)U'(r_0)\geq 0 
\]
holds for $r\in [0,1]$ and $r_0\in (0,1]$.
In addition, for $r,r_0\in (0,1]$, $d_U(r,r_0)=0$ if and only if $r=r_0$.
Recall the limiting behavior of $U$.
\begin{lemma}\label{behaviorU}
The limit
\[
U'(0):=\lim_{h\downarrow0}U'(h)
\]
exists in $[-\infty,\infty)$ and $\lim_{h\downarrow0}h U'(h)=0$ holds.
\end{lemma}
See Appendix~\ref{app:proof4behaviorU} for the proof.
Hereafter, the proofs of the subsequent lemmas are given in Appendix.

By Lemma~\ref{behaviorU}, the limit 
\[
d_U(r,0):=\lim_{r_0 \downarrow0} d_U(r,r_0)\in [0,\infty]
\]
exists.
In the above relation and throughout, we adhere to the following natural convention:
\begin{align*}
u\pm (-\infty)=\mp\infty, \qquad
\lambda \cdot (-\infty)=-\infty, \qquad
-\infty \leq -\infty < u <\infty \leq \infty
\end{align*}
and so on for $u \in \mathbb{R}$ and $\lambda>0$.
Thus, we can regard $d_U$ (resp.\,$D_U$) as a function on $[0,1]\times [0,1]$ (resp. $\mathcal{P}_K \times \mathcal{P}_K$)
valued in $[0,\infty]$.
For $r\in [0,1]$, $d_U(r,0)=0$ if and only if $r=0$.
Moreover, $d_U(r,0)=\infty$ for some $r\in (0,1]$ is equivalent to $U'(0)=-\infty$.

To consider the finiteness of $\mathfrak{D}_U(x,y)$, 
we define the \emph{support} of $z\in \mathcal{P}_K$ by
\[
\spt(z):=\{ k\ |\ z_k>0 \}.
\]
\begin{lemma}\label{support}
For $\Pi\in \Pi(x,y)$, it holds that
\[
    \spt(\Pi)\subset \spt(x)\times \spt(y), \quad \spt(x)\times \spt(y)=\spt(x\otimes y).
\]
\end{lemma}

By Lemma~\ref{support}, we find that $D_U(\cdot, x\otimes y)$ is continuous on a compact set~$\Pi(x,y)$ so that $\mathfrak{D}_U(x,y)<\infty$.
\subsection{On Assumption~\ref{xy} 
with the positivity and finiteness  of $\mathfrak{D}_U(x,y), \Delta_C(x,y)$} \label{subxy}
The problem~\eqref{eq:LP} is trivially solved if 
$\Pi(x,y)=\argmin_{\Pi(x,y)}\langle C,\Pi\rangle$.
Thus, we suppose Assumption~\ref{xy}, in which $\Pi(x,y)$ contains an element other than $x\otimes y$ and hence $\mathfrak{D}_U(x,y)>0$ holds.

Let $V(x,y)$ be the set of vertices of $\Pi(x,y)$.
Under
Assumption~\ref{xy}, 
the set $V(x,y) \setminus \argmin_{V\in V(x,y)}\langle C, V \rangle$ is not empty and $\Delta_C(x,y)\in (0,\infty)$ holds.

\subsection{On Definition~\ref{def:rU}} \label{welldefined}
Let $U$ satisfy Assumption~\ref{assumption}.
It follows from the strict convexity of $U$ on~$[0,1]$ and $U'(0)=-\infty$ that
the function $U':(0,1]\to U'((0,1])=(-\infty, U'(1)]$ admits the inverse function~$e_U$.
We observe from $U''>0$ on $(0,1)$ that 
the function $r\mapsto U'(r)-U'(1-r)$ is strictly increasing on $(0,1)$.
This with the properties 
\[
U'\left(\frac12\right)-U'\left(1-\frac12\right)=0,\qquad
\lim_{r\uparrow 1} ( U'(r)-U'(1-r) )=\infty
\]
guarantees the unique existence of $R_U(x,y)$.
Moreover, since $r\mapsto rU''(r)$ is non-decreasing in $(0,1)$,
we find 
\begin{equation}\label{ru}
\sup_{r \in (0,R_U(x,y)]} \left(U'(1-r)+r U''(r)\right)
\leq U'(1)+R_U(x,y) U''(R_U(x,y))<\infty.
\end{equation}
Thus, all the notions in Definition~\ref{def:rU} are well-defined under Assumption~\ref{assumption}.

\subsection{On Assumption~\ref{assumption}}\label{verify}
By Aleksandrov's theorem (e.g., \cite{EG}*{Theorem~6.9}), 
$U$ is twice differentiable almost everywhere on $[0,1]$.
If $U\in C^2((0,1))$, the strict convexity leads to $U''>0$ almost everywhere on $(0,1)$.
Thus, the requirement $U\in C^2((0,1))$ together with $U''>0$ on $(0,1)$ is mild.

Let $U\in C([0,1])\cap C^1((0,1])\cap C^2((0,1))$ such that $U''>0$ on $(0,1)$.
To apply some algorithms, such as gradient descent~\cite{PeyreCuturi2019}*{Sections~4.4, 4.5, 9.3}, 
we require that 
$\Pi^{U}(\omega)$ belongs to the relative interior of the convex polytope~$\Pi(x,y)$
for any $\omega=(C,x,y,\varepsilon)\in \Omega$,
which holds 
if and only if $U'(0)=-\infty$
by \cite{Tak21arXiv}*{Lemma~3.7 and Remark~3.9}.

Let $U\in C^2((0,1))$ satisfy $U''>0$ on $(0,1)$.
Define $q_U:(0,1)\to [-\infty,\infty]$ by 
\[
q_U(r):=
r U''(r) \cdot \limsup_{h \downarrow 0}
 \frac{1}{h}
 \left(\frac{1}{U''(r+h)}-\frac{1}{U''(r)}\right),
 \qquad
Q_U:=\sup_{r\in (0,1)} q_U(r). 
\]
If $Q_U<\infty$, 
then $U'(0)=-\infty$ yields $Q_U\geq 1$ by \cite{IST}*{Corollaries~2.6, 2.7}.
Note that if $U\in C^3((0,1))$, then
\[
q_U(r)=-\frac{rU'''(r)}{U''(r)}
\quad \text{for }r\in (0,1).
\]
In \cite{IST}, 
the notion of $q_U$ is introduced to determine the hierarchy of $U$ in terms of concavity associated with~$U'$.
See also \cite{OT13}, where $q_U$ is used to classify convex functions into displacement convex classes.
For the definition of the displacement convex classes, 
see~\cite{Vil09}*{Chapter~17}.
It follows from \cite{IST}*{Theorem~2.4} that, 
for $W\in C^2((0,1))$ satisfying $W''>0$ on $(0,1)$,
if $q_U<\infty, q_W>-\infty$ hold 
almost everywhere on $(0,1)$ 
and $q_U \leq q_W$ holds on~$(0,1)$,
then there exist $\lambda>0$ and $\mu \in \mathbb{R}$ such that $U' \geq \lambda W' + \mu$ holds on~$(0,1]$,
consequently, 
\[
e_U(\tau) \leq e_{W}(\lambda^{-1}(\tau-\mu) ) \quad \text{on } \tau \in (U'(0),\lambda W'(1)+\mu.
\]
Thus, under the assumption 
$U\in C([0,1])\cap C^1((0,1])\cap C^2((0,1))$ such that $U''>0$ on ~$(0,1)$ and $U'(0)=-\infty$, 
if $Q_U<\infty$, 
then the choice~$Q_U=1$ is the best possible for the estimate in Theorem~\ref{thm:convergence}.
Moreover,
for $U\in C^2((0,1))$ satisfying  $U''>0$ on $(0,1)$, 
$Q_U=1$ is equivalent to that $r\mapsto r U''(r)$ is non-decreasing on~$(0,1)$ 
by~\cite{IST}*{Corollary~2.6}.
Thus, we confirm that Assumption~\ref{assumption} is reasonable.
Furthermore, our setting is applicable to previous regularization techniques.
Indeed, if we extend a function $U$ satisfying Assumption~\ref{assumption} to a function on $\mathbb{R}$ defined by 
\[
    \widetilde{U}(r) \coloneqq 
    \begin{cases}
        \infty &\text{for } r < 0, \\
        U(r) &\text{for } r\in[0,1], \\
        r\log r + \left( U'(1)-1 \right) r + \left( 1+U(1) - U'(1) \right) &\text{for } r>1,
    \end{cases}
\]
then the convex function on $\mathbb{R}^n$ associated with $\widetilde{U}$ satisfies Dessein et al.'s conditions \cite{DPR18}*{(A) Affine constraints in Table~1}.

In addition, if we choose $W=U_o$, then $q_W\equiv 1$ holds on $(0,1)$, and
hence $Q_U=1$ implies the existence of $\lambda>0$ and $\mu \in \mathbb{R}$ such that
\begin{equation}\label{eu}
e_U(\tau)
\leq \exp\left( \lambda^{-1}(\tau-\mu)\right)
\quad\text{for }\tau \in U'((0,1]).
\end{equation}
Consequently, the error estimate in Theorem~\ref{thm:convergence} is faster than the exponential decay.

\subsection{Asymptotic behavior of the error}\label{error}
Fix $\omega=(C,x,y,\varepsilon)\in \Omega $.
Let $\Pi^\ast\in \Pi(x,y)$ be an optimal transport plan.
We observe from the definition of $\Pi^U(\omega)$ that 
\begin{equation}
\label{seqconv}
\begin{split}
 \langle C, \Pi^U(\omega)\rangle
+\varepsilon D_U(\Pi^U(\omega),x\otimes y)
&=
\inf_{\Pi\in \Pi(x,y)}
\left(
\langle C, \Pi\rangle
+\varepsilon D_U(\Pi,x\otimes y)
\right)\\
&\leq
\langle C,\Pi^\ast\rangle+\varepsilon D_U(\Pi^\ast,x\otimes y).
\end{split}
\end{equation}
This with the nonnegativity of $D_U$ yields
\begin{align*}
\langle C,\Pi^U(\omega)\rangle-\langle C, \Pi^\ast\rangle
&\leq 
\varepsilon
D_U(\Pi^\ast, x\otimes y),
\end{align*}
proving \eqref{conv}.
Moreover, the limit 
\[
\Pi^{U}(C,x,y,0)
:=\lim_{\varepsilon \downarrow 0} 
\Pi^{U}(C,x,y,\varepsilon)
\]
exists and satisfies 
\[
 \argmin_{\Pi'\in\argmin_{\Pi\in \Pi(x,y)}\langle C,\Pi\rangle}D_U(\Pi', x \otimes y)=\{\Pi^{U}(C,x,y,0)\}
\]
(see~\cite{Tak21arXiv}*{Theorem~3.11} for instance).
%

%
\section{Proof of Theorem~\ref{thm:convergence}}
 \label{sec:error_estimate}
Before proving Theorem~\ref{thm:convergence},
we consider the normalization of $U$ 
since the correspondence $U\mapsto D_U$ is not injective.
\begin{lemma}\label{ab}
For $\lambda>0$ and $\mu_0,\mu_1 \in \mathbb{R}$, 
define $U_{\lambda,\mu_0,\mu_1}:[0,1]\to \mathbb{R}$ by 
\[
U_{\lambda,\mu_0,\mu_1}(r):=\lambda U(r)+\mu_1 r+\mu_0.
\]
Then, $d_{U_{\lambda,\mu_0,\mu_1}}=\lambda d_U$ holds on~$[0,1]\times [0,1]$,
consequently,
$D_{U_{\lambda,\mu_0,\mu_1}}=\lambda D_U$ on $\mathcal{P}_K\times \mathcal{P}_K$.
If $U$ satisfies Assumption~\ref{assumption}, then so does $U_{\lambda,\mu_0,\mu_1}$ and $q_{U_{\lambda,\mu_0,\mu_1}}=q_U$ holds on $(0,1)$.
\end{lemma}
Since the proof is straightforward, we omit it.
Let $U$ satisfy Assumption~\ref{assumption}. 
For $\mu_0, \mu_1\in \mathbb{R}$
and $\omega=(C,x,y,\varepsilon) \in \Omega$, 
we have 
\begin{align*}
&\Pi^{U_{1,\mu_0,\mu_1}}(\omega)=\Pi^U(\omega), 
&&\mathfrak{D}_{U_{1,\mu_0,\mu_1}}(x,y)=\mathfrak{D}_U(x,y),\\
&R_{U_{1,\mu_0,\mu_1}}(x,y)=R_U(x,y), 
&&\nu_{U_{1,\mu_0,\mu_1}}(x,y)=\nu_U(x,y)+\mu_1, 
\end{align*}
and $e_{U_{1,\mu_0,\mu_1}}(\tau)=e_U(\tau-\mu_1)$ 
for $\tau\in U_{1,\mu_0,\mu_1}'((0,1])$.
Thus, in Theorem~\ref{thm:convergence}, 
we can normalize $U$ to $U(0)=U(1)=0$ as well as the case of $U=U_o$ without loss of generality.
For the reason and the effect of the choice of $\lambda>0$, see~Section~\ref{UtoL}.

Throughout the rest of this section, we suppose Assumptions~\ref{xy} and \ref{assumption} 
together with $U(0)=U(1)=0$.
We prepare two lemmas to prove Theorem~\ref{thm:convergence}.
Similar statements are found in \cite{Wee18COLT}*{Lemmas~6--8}
but the proofs are drastically different.

\begin{lemma}\label{lem:lower_estiamte_V}
For $r,s,t\in [0,1]$ and $r_0\in (0,1]$, 
\begin{align*}
 U\left((1-t)r+ts\right)&\geq (1-t)U(r)+tU(s)+rU(1-t)+sU(t),\\
d_U( (1-t)r+ ts, r_0)
&\geq (1-t)d_U(r, r_0)
+t d_U(s, r_0)+rU(1-t)+sU(t).
\end{align*}
\end{lemma}

Recall that there exists a unique $R\in (1/2,1)$  such that $U'(R)-U'(1-R)=D$ for any $D>0$
(see~Subsection~\ref{welldefined}).

\begin{lemma}
\label{lem:existence_tb}
For $D>0$ and $R\in (1/2,1)$ with $U'(R)-U'(1-R)=D$, 
the function
\[
r\mapsto Dr-U(r)-U(1-r)
\]
is strictly increasing on $[0,R]$.
Moreover, it holds that
\begin{align}
 -U(r)-U(1-r)
 \le -rU'(r) +r \sup_{\rho\in (0,R]}(U'(1-\rho) +\rho U''(\rho))
 \end{align}
for $r\in (0, R]$.
 \end{lemma}

\textit{Proof of Theorem \ref{thm:convergence}}.
Let $\omega=(C,x,y,\varepsilon)\in \Omega$.
As noted in Subsection~\ref{subxy},
 $\Delta_C(x,y), \mathfrak{D}_U(x,y)\in (0,\infty)$ hold.
We also have
\[
U'(1)\leq U'(1)+\lim_{r \downarrow 0} r U''(r) \leq \nu_U(x,y)<\infty
\]
and $\mathfrak{D}_U(x,y)+\nu_U(x,y)-U'(1)\in (0,\infty)$.
Thus, the interval
\[
\left(0,\frac{\Delta_C(x,y) R_U(x,y)}{\mathfrak{D}_U(x,y)}\right] \cap
\left(0, \frac{\Delta_C(x,y)}{\mathfrak{D}_U(x,y)+\nu_U(x,y)-U'(1)}\right]
\]
is well-defined and nonempty.
Let us choose $\varepsilon$ from the interval.
Note that 
\[
\varepsilon\in
\left(0, \frac{\Delta_C(x,y)}{\mathfrak{D}_U(x,y)+\nu_U(x,y)-U'(1)}\right]
\]
is equivalent to 
\[
-\frac{\Delta_C(x,y)}{\varepsilon}+\mathfrak{D}_U(x,y)
+\nu_U(x,y) \in U'((0,1]).
\]

Recall that $V(x,y)$ is the vertex set of $\Pi (x, y)$.
There exists a family $\{t_V\}_{V \in V(x,y)}\subset [0,1]$  such that 
\[
\sum_{V \in V(x,y)} t_{V}=1, \qquad
\Pi^U(\omega)=\sum_{V \in V(x,y)} t_V V.
\]
Set
\[
V_0(x,y):=\argmin_{V\in V(x,y)}\langle C, V \rangle, \quad
t:=1-\sum_{V\in V_0(x,y)} t_{V}.
\]
By Assumption~\ref{xy},  $V_0(x,y)\neq V(x,y)$ holds.
This with the nonoptimality of  $\Pi^{U}(\omega)$ yields $t>0$.
We also set 
\[
\Pi':=
\sum_{V' \in V(x,y)\setminus V_0(x,y) } \frac{t_{V'}}{t} V',
\quad
\Pi^\ast:=
\begin{dcases}
\sum_{V\in V_0(x,y)} \frac{t_V }{1-t} V &t\neq 1,\\
V^\ast  &t=1,
\end{dcases}
\]
where $V^\ast \in V_0(x,y)$ is chosen arbitrarily.
It turns out that 
$\Pi^\ast, \Pi'\in \Pi(x,y)$ and 
\[
\Pi^U(\omega)=(1-t)\Pi^\ast+ t \Pi'.
\]
We find that 
$
\langle C, V\rangle=\inf_{\Pi\in \Pi(x,y)}
\langle C,\Pi\rangle$
for $V\in V_0(x,y)$, 
in particular
\[
\langle C, \Pi^\ast\rangle=\inf_{\Pi\in \Pi(x,y)}
\langle C,\Pi\rangle.
\]
Setting 
\[
r:=\frac{\langle C, \Pi^U(\omega)\rangle-\langle C,\Pi^\ast \rangle}{\Delta_C(x,y)},
\]
we observe from~\eqref{seqconv}  
with the definition of $\mathfrak{D}_U(x,y)$ that
\begin{align}
\begin{split}
\label{basic}
\langle C, \Pi^U(\omega)\rangle
-\langle C,\Pi^\ast \rangle
&\leq \varepsilon 
\left(
D_U(\Pi^\ast,x\otimes y)
-
D_U(\Pi^U(\omega),x\otimes y)
\right)\\
&\leq \varepsilon \mathfrak{D}_U(x,y).
\end{split}
\end{align}
We also find 
\begin{align*}
\langle C, \Pi^U(\omega)\rangle
-\langle C,\Pi^\ast \rangle
&=t \langle C, \Pi'-\Pi^\ast \rangle\\
&\geq
t\left( \inf_{
V' \in V(x,y)\setminus V_0(x,y)}\langle C,V'\rangle
-
\inf_{V\in V_0(x,y)}\langle C,V\rangle\right)\\
&=t\Delta_C(x,y). 
\end{align*}
These with the condition~$\varepsilon\in (0,\Delta_C(x,y) R_U(x,y)/\mathfrak{D}_U(x,y)]$
yield
\[
t \leq r=\frac{\langle C, \Pi^U(\omega)\rangle
-\langle C,\Pi^\ast \rangle}{\Delta_C(x,y)}
\leq \frac{\varepsilon \mathfrak{D}_U(x,y)}{\Delta_C(x,y)}
\le R_U(x,y).
\]
It follows from Lemma~\ref{support} with the second inequality in Lemma \ref{lem:lower_estiamte_V} that
\begin{align*}
&D_U(\Pi^U(\omega), x\otimes y) \\
& =
D_U( (1-t)\Pi^\ast+ t \Pi',x\otimes y)\\
&=
\sum_{(i,j) \in \spt(x\otimes y)}
d_U( (1-t)\pi^\ast_{ij}+ t \pi'_{ij}, x_iy_j)\\
& \geq
\sum_{(i,j) \in \spt(x\otimes y)}
\left[(1-t)d_U( \pi^\ast_{ij}, x_iy_j)
+
td_U( \pi'_{ij}, x_iy_j)
+\pi^\ast_{ij} U(1-t)+\pi'_{ij}U(t)
\right]\\
&=(1-t)D_U(\Pi^\ast, x\otimes y)
+t D_U(\Pi', x\otimes y)+U(1-t)+U(t).
\end{align*}
This and Lemma~\ref{lem:existence_tb} together with 
$t\leq r\le R_U(x,y)$ yield
\begin{align}
&D_U(\Pi^\ast,x\otimes y)-
D_U(\Pi^U(\omega),x\otimes y)\\
&\leq
tD_U(\Pi^\ast, x\otimes y)-t D_U(\Pi', x\otimes y)
- U(t)-U(1-t)\\
&\leq
t\mathfrak{D}_U(x,y)-U(t)-U(1-t)\\
&\leq
r\mathfrak{D}_U(x,y)-U(r)-U(1-r)\\
&\leq
r\mathfrak{D}_U(x,y)
-rU'(r)+r\nu_U(x,y)\\
&=r(\mathfrak{D}_U(x,y)
-U'(r)+\nu_U(x,y)).
\end{align}
Combining this with \eqref{basic},
we find
\begin{align*}
\frac{\langle C, \Pi^U(\omega)\rangle-\langle C,\Pi^\ast \rangle}{\varepsilon}
&\leq
D_U(\Pi^\ast,x\otimes y)
-
D_U(\Pi^U(\omega),x\otimes y)\\
&\leq
r(\mathfrak{D}_U(x,y)
-U'(r)+\nu_U(x,y)),
\end{align*}
which leads to
\begin{align*}
U'(r)
&\leq 
-\frac{\langle C, \Pi^U(\omega)\rangle-\langle C,\Pi^\ast \rangle}{\varepsilon r}
+\mathfrak{D}_U(x,y)+\nu_U(x,y)\\
&=
-\frac{\Delta_C(x,y)}{\varepsilon }
+\mathfrak{D}_U(x,y)+\nu_U(x,y).
\end{align*}
It follows from the monotonicity of $U'$ on $(0,1]$ that
\begin{align*}
\frac{\langle C, \Pi^U(\omega)\rangle-\langle C,\Pi^\ast \rangle}{\Delta_C(x,y)}
=r
&=e_U(U'(r))\\
&\leq
e_U\left(-\frac{\Delta_C(x,y)}{\varepsilon}+\mathfrak{D}_U(x,y)
+\nu_U(x,y)\right),
\end{align*}
that is,
\begin{align*}
&\langle C, \Pi^U(\omega)\rangle
-\inf_{\Pi\in \Pi(x,y)}\langle C,\Pi\rangle\\
&\leq 
\Delta_C(x,y)
\cdot e_U\left(-\frac{\Delta_C(x,y)}{\varepsilon}+\mathfrak{D}_U(x,y)
+\nu_U(x,y)\right)
\end{align*}
as desired.
$\hfill\qed$

\section{Normalization and scaling}\label{normalization}
In this section, we first show that 
the normalization of a strictly convex function does not affect the error estimate in Theorem~\ref{thm:convergence}.
We then consider the effect of scaling of data and the domain of a strictly convex function  on the error estimate.

\subsection{Normalization}\label{UtoL}
Let $U\in C([0,1]) \cap C^1((0,1])\cap C^2((0,1))$ satisfy $U''>0$ on $(0,1)$.
For $\mu_0,\mu_1\in \mathbb{R}$ and $\lambda>0$, 
set
\[
U_{\lambda,\mu_0,\mu_1}
:=\lambda U(r)+\mu_1 r+\mu_0.
\]
Let $(C,x,y,\varepsilon)\in \Omega$.
Then, 
$\Pi^{U_{\lambda,\mu_0,\mu_1}}(C,x,y,\varepsilon)
=\Pi^{U}(C,x,y,\lambda\varepsilon)$ holds and the interval
\[
\left(0,\frac{\Delta_C(x,y) R_{U_{\lambda,\mu_0,\mu_1}}(x,y)}{\mathfrak{D}_{U_{\lambda,\mu_0,\mu_1}}(x,y)}\right] \cap
\left(0, \frac{\Delta_C(x,y)}{\mathfrak{D}_{U_{\lambda,\mu_0,\mu_1}}(x,y)+\nu_{U_{\lambda,\mu_0,\mu_1}}(x,y)-U_{\lambda,\mu_0,\mu_1}'(1)}\right]
\]
contains $\varepsilon$
if and only if the interval
\[
 \left(0,\frac{\Delta_C(x,y) R_U(x,y)}{\mathfrak{D}_U(x,y)}\right] \cap
\left(0, \frac{\Delta_C(x,y)}{\mathfrak{D}_U(x,y)+\nu_U(x,y)-U'(1)}\right]
\]
contains $\lambda \varepsilon$.
In this case, the equality
\begin{align}
&e_{U_{\lambda,\mu_0,\mu_1}}
\left(-\frac{\Delta_C(x,y)}{\varepsilon}+\mathfrak{D}_{U_{\lambda,\mu_0,\mu_1}}(x,y)
+\nu_{U_{\lambda,\mu_0,\mu_1}}(x,y)\right)\\
&=
e_U\left(-\frac{\Delta_C(x,y)}{\lambda\varepsilon}+\mathfrak{D}_U(x,y)
+\nu_U(x,y)\right)
\end{align}
follows.
Thus, 
we can normalize $U$ as 
\[
U_{\lambda,\mu_0,\mu_1}(0)=
U_{\lambda,\mu_0,\mu_1}(1)=0, \quad
U_{\lambda,\mu_0,\mu_1}'(1)=1,
\]
by choosing $\mu_0, \mu_1\in \mathbb{R}$ and $\lambda>0$ as well as $U_o$.

Let $W\in C([0,1]) \cap C^1((0,1])\cap C^2((0,1))$ also satisfy $W''>0$ on $(0,1)$.
Then, the following four conditions are equivalent to each other.
\begin{enumerate}
\renewcommand{\theenumi}{C\arabic{enumi}}
\renewcommand{\labelenumi}{(\theenumi)}
\setlength{\itemindent}{15pt}
\setcounter{enumi}{-1}
\item \label{C0}
There exist $\mu_0,\mu_1\in \mathbb{R}$ and $\lambda>0$ such that $W=U_{\lambda,\mu_0,\mu_1}$ on $(0,1)$.
\item\label{C1}
There exist 
$\mu_1\in \mathbb{R}$ and $\lambda>0$ such that $W'=\lambda U'+\mu_1$ on $(0,1)$.
\item\label{C2}
There exists $\lambda>0$ such that $W''=\lambda U''$ on $(0,1)$.
\renewcommand{\theenumi}{D}
\renewcommand{\labelenumi}{(\theenumi)}
\item\label{D}
There exists $\lambda>0$ such that $D_W=\lambda D_U$ on $\mathcal{P}_K \times \mathcal{P}_K$ for $K\geq 3$.
\end{enumerate}
We can easily see the equivalence among \eqref{C0}--\eqref{C2}.
The implication from \eqref{C0} to \eqref{D} is straightforward.
Then, here, we only prove that \eqref{D} implies \eqref{C2}.
Assume \eqref{D}.
For $K\geq 3$ and $r\in [0,1]$, define $z^r \in \mathcal{P}_K$ by 
\[
z^r_k:=\begin{cases} 
r & \text{for }k=1,\\
1-r & \text{for }k=2,\\
0& \text{otherwise.}
\end{cases}
\]
We also define $z^\ast\in \mathcal{P}_K$ by $z^\ast_k= K^{-1}$ for all $k$.
For $r\in (0,1)$, we calculate 
\[
W'(r)-W'(1-r)
=\frac{\mathrm{d}}{\mathrm{d}r} D_W(z^r,z^\ast)
=\lambda \frac{\mathrm{d}}{\mathrm{d}r} D_U(z^r,z^\ast)
=\lambda \left( U'(r)-U'(1-r) \right).
\]
Differentiating this with respect to $r$ implies that
\begin{align}\label{UW}
W''(r)-\lambda U''(r)
=
-\left( W''(1-r)-\lambda U''(1-r)\right),
\end{align}
in particular $W''(1/2)=\lambda U''(1/2)$.
For $r,s \in (0,1)$ with $r+s<1$, we define $z^{r,s} \in \mathcal{P}_K$ by 
\[
z^{r,s}_k:=\begin{cases} 
1-(r+s) & \text{for }k=1,\\
r & \text{for }k=2,\\
(K-2)^{-1}s& \text{otherwise.}
\end{cases}
\]
It turns out that 
\begin{align*}
(r+s)W''(1-(r+s))+rW''(r)
&=\frac{\partial}{\partial r} D_W(z^{1},z^{r,s})
=\lambda \frac{\partial}{\partial r} D_U(z^{1},z^{r,s})\\
&=\lambda 
\left[ (r+s)U''(1-(r+s))+rU''(r)\right],
\end{align*}
implying
\begin{align*}
W''(r)-\lambda U''(r)
=-\frac{r+s}{r}\left( W''(1-(r+s))-\lambda U''(1-(r+s))\right).
\end{align*}
This with \eqref{UW} provides 
\[
\frac{r+s}{r}\left( W''(1-(r+s))-\lambda U''(1-(r+s))\right)
=W''(1-r)-\lambda U''(1-r),
\]
in particular, 
the choice of $r=1/2$ leads to
\[
W''\left(\frac12-s\right)=\lambda U''\left(\frac12-s\right) \quad \text{for }s\in \left(0,\frac12\right).
\]
This together with \eqref{UW} gives $W''=\lambda U''$ on $(0,1)$, which is nothing but \eqref{C2}.
Note that, under Assumption~\ref{xy}, we have $I, J\neq 1$ hence $IJ\geq 4$.
Thus, the condition~$K\geq 3$ in \eqref{D} is reasonable. 

We also notice that \eqref{C2} leads to the following condition.
\begin{enumerate}
\setlength{\itemindent}{15pt}
\renewcommand{\theenumi}{C}
\renewcommand{\labelenumi}{(\theenumi)}
\item\label{C4}
$q_U=q_W$ on $(0,1)$.
\end{enumerate}
By \cite{IST}*{Theorem~2.4}, if $q_U,q_W$ are finite almost everywhere on $(0,1)$,
then \eqref{C4} leads to \eqref{C2}. 
Thus, all conditions~\eqref{C0}--\eqref{C2}, \eqref{D}, and \eqref{C4} are equivalent to each other.

By the equivalence between \eqref{C0} and \eqref{C1},
choosing either $U$ or $U'$ determines the other under the normalization
\begin{equation}\label{001}
U(0)=U(1)=0, \qquad U'(1)=1.
\end{equation}
In what follows, we use the symbol~$L$ instead of $U'$ and interpret Assumption~\ref{assumption} for $U$ as follows.
\begin{assumption}\label{ell}
Let $L\in C((0,1])\cap C^1((0,1))$ satisfy $L'>0$ on $(0,1)$ and  $\lim_{t\downarrow 0}L(t)=-\infty$.
In addition, $t \mapsto tL'(t)$ is non-decreasing on $(0,1)$.
\end{assumption}
Suppose Assumption~\ref{ell}.
Then, there exists $t_0\in (0,1]$ such that $L<0$ on $(0,t_0]$.
By the monotonicity of $t\mapsto tL'(t)$, we have
 \[
L(t_0)-L(t)=\int_t^{t_0} L '(s) \mathrm{d}s\leq t_0L'(t_0) \int_{t}^{t_0} \frac{1}{s} \mathrm{d}s=t_0L'(t_0) \left(\log t_0-\log t\right)
 \]
for $t\in (0,t_0]$.
Since $L$ is monotone increasing on $(0,t_0]$ and 
 \[
\left( L(t_0)-t_0 L'(t_0) \log t_0 \right) (t_0-h)+ t_0L'(t_0)\int_h^{t_0} \log t\mathrm{d}t
\leq \int_h^{t_0}L(t)\mathrm{d}t <0
 \]
holds $h \in (0,t_0]$,
the improper integral
\[
U_L(r):=\int_0^r L(t)\mathrm{d}t
\]
is well-defined for $r\in [0,1]$. 
It is easy to see that $U_L$ satisfies Assumption~\ref{assumption}.
Note that 
\[
d_{U_L}(r,r_0)
=\int_{r_0}^r (L(t)-L(r_0)) \mathrm{d}t
=\int_{r_0}^r \int_{r_0}^t L'(s)  \mathrm{d}s\mathrm{d}t
\quad \text{for }r\in [0,1], r_0\in (0,1].
\]
Conversely, if $U$ satisfies Assumption~\ref{assumption},
then $L=U'$ satisfies Assumption~\ref{ell}.
Thus, we can use $L$ satisfying Assumption~\ref{ell}
instead of $U$ satisfying Assumption~\ref{assumption} for our regularization.

\begin{remark}\label{remark4}
In Theorem~\ref{thm:convergence}, 
the range of the regularization parameter~$\varepsilon$ is given by the intersection of the two intervals.
One interval
\[
\left(0, \frac{\Delta_C(x,y)}{\mathfrak{D}_U(x,y)+\nu_U(x,y)-U'(1)}\right]
\]
is needed to make
\[
-\frac{\Delta_C(x,y)}{\varepsilon}+\mathfrak{D}_U(x,y)
+\nu_U(x,y)\in U'((0,1])
\]
as seen in the proof of Theorem~\ref{thm:convergence}.
Hence, this interval is not needed if the function~$U$ is extended to a continuous, strictly convex function on $[0,\infty)$ and $U\in C^1((0,\infty))$ with $\lim_{r\uparrow \infty} U'(r)=\infty$.
If $U$ satisfies the the normalization~\eqref{001},
then
we can extend $U$ to as a function on $[0,\infty)$ by setting 
\[
U(r):=r \log r \quad \text{for }r>1.
\]

For $L$ satisfying Assumption~\ref{ell}, 
if we set 
\[
\ell(t):=\left(L(1)-\int_{0}^1 L(s) \mathrm{d}s \right)^{-1} \left(L(t)-\int_{0}^1 L(s) \mathrm{d}s\right)
\]
for $t\in [0,1]$, then $\ell$ satisfies Assumption~\ref{ell} and $U_\ell$ satisfies the normalization~\eqref{001}.
\end{remark}

\begin{remark}
We must take care of the normalization of $C\in \mathbb{R}^{I\times J}$, because the normalization on $C$ and the scale of the regularization parameter~$\varepsilon>0$ are in the relationship of trade-off, that is,
\[
\Pi^U(C,x,y,\varepsilon) 
=
\Pi^U(\lambda C,x,y,\lambda^{-1}\varepsilon).
\]
\end{remark}

\subsection{Scaling}\label{well-defined}
In Weed's work~\cite{Wee18COLT}, 
the components of vectors $x,y$ are only assumed to be nonnegative.
To extend our estimate to such a general setting,
let us consider scaling vectors 
and show that this is equivalent to scaling the domain of a strictly convex function.

In the optimal transport problem~\eqref{eq:LP}, 
although two given data~$x$ and $y$ are normalized to be 1 with respect to the $\ell^1$-norm,
their $\ell^1$-norms can be chosen arbitrarily if both are the same.
For a subset~$\mathcal{Z}$ of Euclidean space and $a>0$,
set 
\[
a \mathcal{Z}:=\{az\bigm| z\in \mathcal{Z} \}.
\]
For $x\in \mathcal{P}_I$ and $y\in \mathcal{P}_J$, 
we shall, by abuse
of notation, define 
\begin{align*}
 &\Pi(ax,ay)\\&:=
 \left\{
\widetilde{\Pi}=(\widetilde{\pi}_{ij})\in a\mathcal{P}_{I \times J}
 \biggm|
 \sum_{l=1}^{J}\widetilde{\pi}_{il}=ax_i \text{ and } 
 \sum_{l=1}^{I}\widetilde{\pi}_{lj}=ay_j\ 
 \text{ for any } i,j
 \right\}.
\end{align*}
Then, we have $\Pi(ax, ay)=a \Pi(x,y)$ and hence $a x\otimes y\in \Pi(ax,ay)$.
We denote by $V(ax,ay)$ the set of the vertices of $\Pi(ax,ay)$. 
Then, $V(ax,ay)=aV(x,y)$ follows.

For $U\in C([0,1])\cap C^1((0,1])\cap C^2((0,1))$ such that $U''>0$ on $(0,1)$
and $a'\in (0,1]$, 
we can define $D_{U} :a'\mathcal{P}_K \times a'\mathcal{P}_K \to [0,\infty]$ by
\[
D_{U}(a'z,a'w):=\sum_{k} d_{U}(a'z_k, a'w_k)
\quad\text{for } z,w\in \mathcal{P}_K
\]
and consider the regularized problem
\begin{equation}\label{eq:Otee}
\inf_{\widetilde{\Pi}\in \Pi(a'x,a'y)} 
\left(
\langle C,\widetilde{\Pi}\rangle+\varepsilon D_{U}(\widetilde{\Pi}, a'x\otimes y)
\right)
\quad \text{for }\omega=(C,x,y,\varepsilon)\in \Omega.
\end{equation}
More generally, for $W\in C([0,b])\cap C^1((0,b])
\cap C^2((0,b))$ with $b>0$ such that $W''>0$ on $(0,b)$,
we define $d_W:[0,b] \times [0,b] \to [0,\infty]$ by
\begin{align*}
 d_W(r,r_0)
:= W(r)-W(r_0)-(r-r_0)W{}'(r_0) \quad\text{for }r\in [0,b], r_0\in (0,b]
\end{align*}
and 
\[
 d_W(r, 0):=\lim_{h \downarrow 0}d_W(r,h)
\quad\text{for }r\in [0,b].
\]
For $0<b'\leq b$,
we  define 
$D_{W} :b'\mathcal{P}_K \times b'\mathcal{P}_K \to [0,\infty]$ by
\begin{align*}
D_{W}(b'z,b'w)&:=\sum_{k} d_{W}(b'z_k, b'w_k)
\quad\text{for } z,w\in \mathcal{P}_K.
\end{align*}
This enables us to consider a similar regularized problem on $b'\mathcal{P}_I \times b'\mathcal{P}_J$ to~\eqref{eq:Otee}.
Also, we can scale the domain of $U$. 
For $a>0$, define $W_b^a$ by
\begin{equation}\label{wab}
W_b^a(r):=ab^{-1}W(a^{-1} br)
\quad \text{for }r\in [0,a]. 
\end{equation}
We find $W_b^a \in C([0,a]) \cap C^1((0,a])\cap C^2((0,a))$ with $W_b^a{}''>0$ on $(0,a)$.
Moreover,
if 
$\lim_{h\downarrow0}W'(h)=-\infty$ holds 
and $r \mapsto r W''(r)$ is non-decreasing on~$(0,b)$,
then 
$\lim_{h\downarrow0}W_b^a{}'(h)=-\infty$ holds
and $r \mapsto r W_b^a{}''(r)$ is non-decreasing on~$(0,a)$.
In particular, $U:=W_b^1$ satisfies Assumption~\ref{assumption}.
Note that the condition
$W(0)=W(b)=0$ with $W'(b)=1$
is equivalent to the condition
$W_b^a(0)=W_b^a(a)=0$ with $W_b^a{}'(a)=1$.
We observe from 
\begin{align*}
d_{W^b_a}(r,r_0)=a\cdot d_U(a^{-1}r, a^{-1}r_0)
\quad
\text{for } r,r_0\in [0,a]
\end{align*}
that
\begin{align}
&\inf_{\widetilde{\Pi}\in \Pi(ax,ay)} 
\left(
\langle C,\widetilde{\Pi}\rangle+\varepsilon D_{W^a_b}(\widetilde{\Pi}, ax\otimes y)
\right)\\
&=
a\cdot \inf_{\Pi\in \Pi(x,y)} 
\left(
\langle C, \Pi\rangle+ \varepsilon D_{U}(\Pi, x\otimes y)
\right),\\
\begin{split}
&\argmin_{\widetilde{\Pi}\in \Pi(ax,ay)} 
\left(
\langle C,\widetilde{\Pi}\rangle+\varepsilon D_{W^a_b}(\widetilde{\Pi}, ax\otimes y)
\right)\\
&=
a \cdot \argmin_{\widetilde{\Pi}\in \Pi(x,y)} 
\left(
\langle C,\widetilde{\Pi}\rangle+ \varepsilon D_{U}(\widetilde{\Pi}, a^{-1} x\otimes y)
\right)
\end{split} \label{eq:argmin}
\end{align}
for $\omega=(C,x,y,\varepsilon)\in \Omega$.
This means that the scaling of $U$ cancels out the scaling of $x$ and $y$.
Further, we see that the set~\eqref{eq:argmin} is a singleton 
and we denote by $\Pi^{W_b^a}(C,ax,ay,\varepsilon)$ the unique element.
Then, $\Pi^{W_b^a}(C,ax,ay,\varepsilon)=a\Pi^U(\omega)$ holds. 
Let us define all notions required to state Theorem~\ref{thm:convergence} 
as follows:
\begin{equation}
\label{eq:quantity}
\begin{split}
 \mathfrak{D}_{W_b^a}(ax, ay)
&\coloneqq\sup_{\widetilde{\Pi}\in\Pi(ax, ay)}D_{W_b^a}(\widetilde{\Pi},a x\otimes y),\\
\Delta_C(ax, ay)
 &\coloneqq
 \inf_{\widetilde{V}' \in V(ax,ay) \setminus \argmin_{\widetilde{V}\in V(ax, ay)} \langle C,\widetilde{V}'\rangle }\langle C,\widetilde{V}'\rangle -\inf_{\widetilde{V}\in V(ax, ay)} \langle C, \widetilde{V} \rangle,\\
R_{W_b^a}(ax, ay)&\in \left[a/2,a\right]\ \text{ such that }\\ 
&W_b^a{}'(R_{W_b^a}(ax, ay))-W_b^a{}'(a-R_{W_b^a}(ax, ay))=a^{-1}\mathfrak{D}_{W_b^a}(ax, ay),\\
\nu_{W_b^a}(ax, ay)&:=\sup_{r \in (0,R_{W_b^a}(ax, ay)]} \left(W_b^a{}'(a-r)+r W_b^a{}''(r)\right).
\end{split}
\end{equation}
Let $e_{W_b^a}$ be the inverse function of $W_b^a{}':(0,a]\to W_b^a{}'((0,a])$.
%
It follows from 
\begin{alignat}{2}
W_b^a(ar)
&=aU(r), 
&\qquad
W_b^a{}'(ar)
&=U'(r), \\
W_b^a{}''(ar)
&=a^{-1}U''(r), 
&\qquad
d_{W_b^a}(ar,ar_0)
&=a d_U(r,r_0)
\end{alignat}
for 
$r,r_0\in [0,1]$ that
\begin{align*}
%
\mathfrak{D}_{W_b^a}(ax,ay)&=a\mathfrak{D}_U(x,y), 
&\qquad
\Delta_C(ax, ab)&=a \Delta_C(x,y),\\
R_{W_b^a}(ax,ay)&=aR_U(x,y), 
&\qquad
\nu_{W_b^a}(ax,ay)&=\nu_U(x,y),
\end{align*}
and $e_{W_b^a}=ae_U$  on  $W_b^a{}'((0,a])=U'((0,1])$.
Thus, under Assumption~\ref{xy}, the interval 
\begin{align*}
& \left(0,\frac{\Delta_C(ax,ay) R_{W_b^a}(ax,ay)}{a\mathfrak{D}_{W_b^a}(ax,ay)}\right]\\ 
&\cap
\left(0, \frac{a^{-1}\Delta_C(ax,ay)}{a^{-1}\mathfrak{D}_{W_b^a}(ax,ay)+\nu_{W_b^a}(ax,ay)-W_b^a{}'(a)}\right]
\end{align*}
is determined independently of the choice of $a,b$, 
which coincides with 
\begin{align*}
\left(0,\frac{\Delta_C(x,y) R_U(x,y)}{\mathfrak{D}_U(x,y)}\right] \cap
\left(0, \frac{\Delta_C(x,y)}{\mathfrak{D}_U(x,y)+\nu_U(x,y)-U'(1)}\right].
\end{align*}
In addition, we also have
\begin{align*}
&\langle C, \Pi^{W_b^a}(C,ax,ay,\varepsilon)\rangle
-\inf_{\widetilde{\Pi} \in \Pi(ax,ay)}\langle C,\widetilde{\Pi}\rangle\\
&=
\langle C, a\Pi^U(\omega)\rangle
-\inf_{\Pi \in \Pi(x,y)}\langle C,a\Pi\rangle\\
&=
a\left( 
\langle C, \Pi^{U}(\omega)\rangle
-\inf_{\Pi \in \Pi(x,y)}\langle C,\Pi\rangle\right)\\
&\leq 
a\left( 
\Delta_C(x,y)
e_U\left(-\frac{\Delta_C(x,y)}{\varepsilon}+\mathfrak{D}_U(x,y)
+\nu_U(x,y)\right)\right)\\
&=
\Delta_C(x,y)
e_{W_b^a}\left(-\frac{\Delta_C(x,y)}{\varepsilon}+\mathfrak{D}_U(x,y)
+\nu_U(x,y)\right)\\
&=
\frac{\Delta_C(ax,ay)}{a}
e_{W_b^a}\left(-\frac{\Delta_C(ax,ay)}{a\varepsilon}+\frac1a\mathfrak{D}_{W_b^a}(ax,ay)
+\nu_{W_b^a}(ax,ay)\right)
\end{align*}
for $\varepsilon$ in the interval above.
This estimate can be derived directly in a similar way to the proof of Theorem~\ref{thm:convergence}.

Thus, in the problem~\eqref{eq:OTe},
if we simultaneously scale the data and the domain of a strictly convex function,
we obtain essentially the same error estimate in Theorem~\ref{thm:convergence}, where the domain of a strictly convex function has no effect on the estimate.

\subsection{Invariance of the Kullback--Leibler divergence under scaling of data}
\label{scalenormalized}
In contrast to the previous subsection, let us now consider the case where
the scaling of data does not match the domain of a strictly convex function.
For this sake, 
fix $\widetilde{a},a,b>0$ with $\widetilde{a}<a$.
Let 
$W\in C([0,b]) \cap C^1((0,b])\cap C^2((0,b))$ satisfy that
$W''>0$ on $(0,b)$, $\lim_{h\downarrow0}W'(h)=-\infty$, 
and $r \mapsto r W''(r)$ is non-decreasing on $(0,b)$ as in the previous subsection. 

First, we find that
\begin{align*}
&\argmin_{\widetilde{\Pi}\in \Pi(\widetilde{a}x,\widetilde{a}y)}
\left( \langle C,\widetilde{\Pi} \rangle+ \varepsilon D_{W_b^a}(\widetilde{\Pi}, \widetilde{a}x \otimes y)\right)\\
&=
\argmin_{\widetilde{\Pi}\in \Pi(\widetilde{a}x,\widetilde{a}y)}
\left( \langle C,a^{-1}\widetilde{\Pi} \rangle+ \varepsilon D_{W_b^{1}} (a^{-1}\widetilde{\Pi}, a^{-1}\widetilde{a} x \otimes y)\right)\\
&=
\argmin_{\Pi\in \Pi(x,y)}
\left( \langle C,\Pi \rangle+ \varepsilon D_{W_b^{a/ \widetilde{a}}} (\Pi, x \otimes y)\right)
\end{align*}
for $\omega=(C,x,y,\varepsilon)\in\Omega$.
Thus, scaling the data is equivalent to scaling the domain of a strictly convex function. 
Then, we can use $U$ satisfying Assumption~\ref{assumption} instead of $W$.

The following proposition suggests that the regularization effect by a Bregman divergence varies under scaling of data unless the Bregman divergence is the Kullback--Leibler divergence.

\begin{proposition}
\label{prop:DU}
Let $U$ satisfy Assumption~\ref{assumption}.
Let $\widetilde{a},a>0$ with $\widetilde{a}<a\leq 1$.
If there exists 
$\kappa>0$ such that 
\begin{equation}\label{K}
D_U(\widetilde{a}z,\widetilde{a}w)=\kappa D_U(az,aw) \quad\text{for } z,w\in \mathcal{P}_K \quad \text{for } K\geq 3,
\end{equation}
then there exist $\mu_0,\mu_1\in \mathbb{R}$ and $\lambda>0$ such that $U=(U_o)_{\lambda,\mu_0,\mu_1}$ on $(0,a]$.
\end{proposition}

Let $U$ satisfy Assumption~\ref{assumption} 
and $a \in (0,1)$.
For a regularized problem of the form
\[
\inf_{\widetilde{\Pi}\in \Pi(\widetilde{a}x,\widetilde{a}y)}
\left( \langle C,\widetilde{\Pi} \rangle+ \varepsilon D_{U} (\widetilde{\Pi}, a x \otimes y)\right)
\quad
\text{for } \omega=(C,x,y,\varepsilon)\in \Omega,
\]
the choice of $a$ is important 
to give a similar estimate as in Theorem~\ref{thm:convergence},
since the quantities such as \eqref{eq:quantity}
may be involved in the estimate.
Indeed, 
if we define
\[
 \mathfrak{D}_{U}(ax, ay)
\coloneqq\sup_{\widetilde{\Pi}\in\Pi(ax, ay)}D_{U}(\widetilde{\Pi},a x\otimes y)
\]
then, for $\widetilde{a}>0$,
it turns out that 
\begin{align*}
d_{U}(ar,ar_0)
&=ad_{U^{a^{-1}}}(r,r_0)\\
&=a \int_{r_0}^{r} \int_{r_0}^t  \frac{\mathrm{d}^2}{\mathrm{d}s^2}U^{a^{-1}}{}(s) \mathrm{d}s\mathrm{d}t 
=a\int_{r_0}^r \int_{r_0}^t  a U''(a s) \mathrm{d}s \mathrm{d}t \\
&\geq
a
\int_{r_0}^r \int_{r_0}^t \widetilde{a}  U''(\widetilde{a} s) \mathrm{d}s \mathrm{d}t 
=a\widetilde{a}^{-1}
d_{U}(\widetilde{a}r,\widetilde{a}r_0),
\end{align*}
consequently, 
\[
\widetilde{a}^{-1}\mathfrak{D}_{U}(\widetilde{a}x,\widetilde{a}y) \leq 
a^{-1} \mathfrak{D}_{U}(ax,ay)
\]
holds with equality if and only if 
$U=(U_o)_{\lambda,\mu_0,\mu_1}$
holds for some $\mu_0, \mu_1\in \mathbb{R}$ and $\lambda>0$.
\section{Examples and comparison}\label{example} 
We give examples of $U$ satisfying Assumption~\ref{assumption}.
\subsection{Model case}\label{model}
Recall that $U_o\in C([0,\infty))\cap C^{\infty}((0,\infty))$ is defined as 
\[
U_o(r):=
\begin{cases}
r \log r &\text{for }r>0,\\
0& \text{for }r=0.
\end{cases}
\]
Obviously, $U_o$ satisfies Assumption~\ref{assumption} and the normalization~\eqref{001}.
For $r,r_0>0$, we see that 
\[
d_{U_o}(r,r_0)
=
U_o(r)-U_o(r_0)-(r-r_0)U_o'(r_0)
=r (\log r-\log r_0)-(r-r_0),
\]
which yields $D_{U_o}(z,w)=D_{\mathrm{KL}}(z,w)$ for $z, w\in \mathcal{P}_K$.

Let us see that Theorem~\ref{thm:convergence} for the case $U=U_o$ coincides with the error estimate given in \cite{Wee18COLT}*{Theorem~5}
interpreted as 
 \begin{align}\label{weed}
\langle C, \Pi^{U_o}(\omega)\rangle-\inf_{\Pi\in \Pi(x,y)}\langle C,\Pi\rangle	 \leq \Delta_C(x,y) \exp \left( - \frac{\Delta_C(x,y)}{\varepsilon} + \mathfrak{D}_{U_o}(x,y) + 1 \right)
\end{align}
for $\varepsilon \in \left(0,{\Delta_C(x,y)}/(1+\mathfrak{D}_{U_o}(x,y))\right]$,
where $\omega=(C,x,y,\varepsilon)\in \Omega$.
In our setting, 
the $\ell^1$-radius of $\Pi(x,y)$ defined in \cite{Wee18COLT}*{Definition~2} is calculated as 
\[
\max_{\Pi \in \Pi(x,y)} \sum_{i,j} \pi_{ij}=1.
\]
The entropic radius defined in \cite{Wee18COLT}*{Definition~3} is calculated as
\begin{align*}
&\sup_{\Pi,\Pi'\in\Pi(x,y)} \left(S(\Pi) - S(\Pi')\right) \\
& = \sup_{\Pi,\Pi'\in\Pi(x,y)}\left( -D_{U_o}(\Pi, x\otimes y)+D_{U_o}(\Pi', x\otimes y)\right)\\
& = \sup_{\Pi\in\Pi(x,y)}D_{U_o}(\Pi, x\otimes y)-\inf_{\Pi\in\Pi(x,y)}D_{U_o}(\Pi, x\otimes y) \\
&=\mathfrak{D}_{U_o}(x,y),
 \end{align*}
thanks to the relation 
$ D_{U_o}(\Pi,x\otimes y)=-S(\Pi)+S(x)+S(y)$ for $\Pi \in \Pi(x,y)$.
A direct calculation provides 
\begin{alignat*}{2}
e_{U_o}(\tau)&=\exp(\tau-1):(-\infty,1]\to (0,1], 
&\qquad R_{U_o}(x,y)&=\frac{e^{\mathfrak{D}_{U_o}(x,y)}}{1+e^{\mathfrak{D}_{U_o}(x,y)}},\\
\nu_{U_o}(x,y)&=2, 
&\qquad U_o'(1)&=1.
\end{alignat*}
Thus, our estimate in Theorem~\ref{thm:convergence} coincides with \eqref{weed},
where the range of the regularization parameter in Theorem~\ref{thm:convergence} is given by 
\begin{align*}
& \left(0,\frac{\Delta_C(x,y) R_{U_o}(x,y)}{\mathfrak{D}_{U_o}(x,y)}\right] \cap
\left(0, \frac{\Delta_C(x,y)}{\mathfrak{D}_{U_o}(x,y)+\nu_{U_o}(x,y)-U_o'(1)}\right]\\
&=
\left(0,\frac{\Delta_C(x,y)}{1+\mathfrak{D}_{U_o}(x,y)}\right],
\end{align*}
which coincides with that in \cite{Wee18COLT}*{Theorem~5}.

Note that a symmetric version of $U_o$ defined by
\begin{equation}
U_{\mathrm{s}}(r):=  U_o(r)+U_o(1-r) \quad \text{for }r\in [0,1]
\label{sym}
\end{equation}
satisfies Assumption~\ref{assumption} except for the differentiability at $r=1$.
The entropy associated with the function $U_{\mathrm{s}}$ is called the \emph{Fermi--Dirac entropy}.
As mentioned in Subsection~\ref{scalenormalized},
we can scale $U_{\mathrm{s}}$ to satisfy Assumption~\ref{assumption}.

\subsection{\texorpdfstring{$q$}{q}-logarithmic function}
Let us consider an applicable example other than the model case.
From the equivalence between \eqref{C0}--\eqref{C2}, any of $U, U'$, and $U''$ is on the table for consideration.
If we regard~$1/U''$ as a deformation function,  $U'$ is called a deformed logarithmic function
and $-U$ corresponds to the density function of an entropy (see~\cite{Naudts}*{Chapters~10, 11} for details).
One typical example of deformed logarithmic functions is the $q$-logarithmic function.
For $q\in \mathbb{R}$, 
define the \emph{$q$-logarithmic function} $\ln_q\colon(0,\infty)\to \mathbb{R}$ by 
 \begin{align}
 \ln_q(t)
 \coloneqq
 \int_1^t s^{-q}\,\mathrm{d}s
 =
 \begin{dcases*}
 \frac{t^{1-q}-1}{1-q}&if $q\neq1$,\\
 \log t&if $q=1$.
 \end{dcases*}
 \end{align}
The entropy associated with $\ln_q$ is called the Tsallis entropy (see~\cite{Naudts}*{Chapter~8} for instance).
We see that the property $\lim_{h \downarrow 0}\ln_q(h)=-\infty$ is equivalent to $q \geq 1$.
Moreover, the function $t \mapsto t \ln_q'(t)=t^{1-q}$ is non-decreasing on~$(0,1)$ if and only if $q\leq 1$ holds.
Thus, $\ln_q$ satisfies Assumption~\ref{ell} if and only if $q=1$ holds,
where $U_{\ln_1}=(U_o)_{1, 0, -1}$.
This means that, 
if we regard a power function as a deformation function, that is, $1/U''$, 
only the power function of exponent~$1$ is applicable to Theorem~\ref{thm:convergence}.
Thus, for some existing studies including the Tsallis entropy, an error estimate similar to the case of the Kullback--Leibler divergence would not be obtained.
To obtain an example that satisfies the conditions of Theorem~\ref{thm:convergence},
we need modify the power function of exponent~$1$ other than power functions of general exponent.

\subsection{Upper incomplete gamma function}\label{powerlog}
For $\alpha\in \mathbb{R}$,
define $L_\alpha:(0,1)\to \mathbb{R}$ by 
\[
L_\alpha(t):=-\ln_{1-\alpha}(-\log t)
=\begin{cases}
-\dfrac{(-\log t)^\alpha-1}\alpha   &\text{if }\alpha \neq 0,\\
-\log (-\log t) &\text{if }\alpha=0.
\end{cases}
\]
It turns out that 
\[
\frac1{L_{\alpha}'(t)}= t\cdot (-\log t)^{1-\alpha}>0
\quad \text{for }t\in (0,1),
\]
which can be regarded as a refinement of the power function of exponent~$1$ 
since the logarithmic function is referred to as the power function of exponent~$0$. 
We see that 
$\lim_{t \uparrow 1}L_\alpha(t)$ is finite 
if and only if $\alpha>0$ holds,
and 
$\lim_{h \downarrow 0}L_\alpha(h)=-\infty$
if and only if $\alpha \geq 0$ holds.
Moreover, the function $t\mapsto t L'_\alpha(t)=(-\log t)^{\alpha-1}$ is non-decreasing on $(0,1)$ 
if and only if $\alpha \leq 1$ holds.
Thus, $L_\alpha$ satisfies Assumption~\ref{ell} if and only if $\alpha\in(0,1]$ holds,
where $U_{L_1}=U_o$.

In what follows, $\alpha\in(0,1]$ is assumed.
We set 
\[
L_\alpha(1):=\lim_{t\uparrow 1} L_\alpha(t)=\frac1{\alpha}.
\]
It follows from 
the change of variables $-\log t = \tau$ that
\[
\int_{0}^1 L_\alpha(t) \mathrm{d}t
=-\frac{1}{\alpha}\Gamma(\alpha+1)+\frac1\alpha,
\]
where $\Gamma(\cdot)$ is the gamma function.
As mentioned in Remark~\ref{remark4},
if we set
\[
\ell_{\alpha}(t)
:=-\frac{(-\log t)^\alpha}{\Gamma(\alpha+1)}+1,
\]
for $t\in (0,1]$, then $\ell_\alpha$ satisfies Assumption~\ref{ell}
and
\begin{equation}\label{test}
U_{\alpha}(r)
:=\int_0^r \ell_\alpha(t)\mathrm{d}t
=-\frac{1}{\Gamma(\alpha+1)}\Gamma(\alpha+1, -\log r) +r
\end{equation}
satisfies the normalization~\eqref{001},
where
\[
\Gamma(p, \tau) = \int_\tau^\infty t^{p-1} \exp(-t) \mathrm{d}t
\quad\text{for }\tau \geq0
\]
is the upper incomplete gamma function for $p > 0$.
Note that $\Gamma(p,0)=\Gamma(p)$.
Since the inverse function 
$e_{U_\alpha}$ of $U_{\alpha}':(0,1]\to U'_\alpha((0,1]) =(-\infty, 1]$ is given by 
\[
e_{U_\alpha}(\tau)=\exp\left(-\left[ -\Gamma(\alpha+1)(\tau-1)\right]^{\frac1\alpha}\right),
\]
the error estimate in Theorem~\ref{thm:convergence} for $U=U_\alpha$ 
decays exponentially when $\alpha=1$, which is the same as~\eqref{weed}, 
and is more tight if $\alpha \in (0,1)$ as  observed in Subsection~\ref{verify}.

The function $L_\alpha$ is introduced to analyze the preservation of concavity by
the Dirichlet heat flow in a convex domain on Euclidean space
(see \cites{IST01,IST07}).
\subsection{Complementary error function} \label{erfc}
Let us give an example of a function $L$ satisfying Assumption~\ref{ell} except for the continuity at $t=1$.
In this case, 
\[
W(r):=\int_0^r L(t)\mathrm{d}t
\]
satisfies Assumption~\ref{assumption} except for the continuity and differentiability at $r=1$
but 
\[
W^a_1(r)
=aW(a^{-1}r)
=a \int_0^{a^{-1}r} L(t)\mathrm{d}t
=\int_0^{r} L(a^{-1}t)\mathrm{d}t
\]
satisfies Assumption~\ref{assumption} if $a>1$.
Then, it might be worth considering the effect of $a>1$ on the regularized solution of the form
\begin{align*}
&\quad\argmin_{\Pi \in \Pi(x, y)} 
\left( \langle C, \Pi \rangle+\varepsilon D_{W_1^a}(\Pi,  x\otimes y)\right),
\end{align*}
as we mentioned in Subsection~\ref{scalenormalized}.

Let us give an example of such $L$.
Define $H:(0,1)\to \mathbb{R}$ by the inverse function of 
\[
\tau\mapsto \frac{1}{\sqrt{\pi}} \int_{-\tau/2}^{\infty} e^{-\sigma^2} d\sigma
=\frac12\operatorname{erfc}\left(-\frac{\tau}2\right),
\]
where $\operatorname{erfc}(\tau) = 1-\operatorname{erf}(\tau)$ is the complementary error function with the error function
\begin{align*}
    \operatorname{erf}(\tau) = \frac{2}{\sqrt{\pi}}\int_0^\tau e^{-\sigma^2} d\sigma
\quad \text{for }\tau \in \mathbb{R},
\end{align*}
and we used the properties 
\[
\lim_{\tau \downarrow-\infty}
\frac12\operatorname{erfc}\left(-\frac{\tau}2\right)=0, \qquad
\lim_{\tau \uparrow \infty}
\frac12\operatorname{erfc}\left(-\frac{\tau}2\right)=1.
\]
It is easily seen  that $H\in C^\infty((0,1))$ and $\lim_{t\downarrow 0} H(t)=-\infty$.
Since we have 
\begin{align*}
\frac12\operatorname{erfc}\left(-\frac{H(t)}{2}\right)&=t \quad \text{for }t\in (0,1), \\
\frac{\mathrm{d}}{\mathrm{d}\tau}\left(\frac12\operatorname{erfc}\left(-\frac{\tau}2\right) \right)&=\frac{1}{\sqrt{4\pi}}e^{-\frac{\tau^2}{4}}
\quad \text{for }\tau \in \mathbb{R},
\end{align*}
we find that 
\begin{align*}
1&=\frac{\mathrm{d}}{\mathrm{d}t}\left( \frac12\operatorname{erfc}\left( -\frac{H(t)}2 \right)\right) =\frac{1}{\sqrt{4\pi}}  e^{-\frac{H(t)^2}{4}}H'(t),\\
0&=\frac{\mathrm{d}^2}{\mathrm{d}t^2} \left( \frac12\operatorname{erfc}\left( -\frac{H(t)}2 \right)\right)=\frac{1}{\sqrt{4\pi}}  e^{-\frac{H(t)^2}{4}}
\left( -\frac{H(t)}{2}H'(t)^2+H''(t) \right),
\end{align*}
consequently,
\begin{align*}
\frac{\mathrm{d}}{\mathrm{d}t}(t H'(t))
&=
H'(t)+tH''(t)
=H'(t)\left(1+\frac{tH(t)}{2} H'(t) \right),
\end{align*}
for $t\in (0,1)$.
It is proved in \cite{IST}*{Section~4.3} that 
\[
\inf_{t\in (0,1)}
\frac{tH(t)}{2}H'(t) 
=
\lim_{t\downarrow 0}
\frac{tH(t)}{2} H'(t)
=-1
\]
in terms of the inverse function of $H$,
and hence $t \mapsto t H'(t)$ is non-decreasing on $(0,1)$.
Thus, $H$ satisfies Assumption~\ref{ell} except for the continuity at $t=1$.

The function $H$ is also introduced to analyze the preservation of concavity by
the Dirichlet heat flow in a convex domain on Euclidean space
(see~\cite{IST07}).
For $F \in C^2((0,1))$, 
if $F$-concavity is preserved by the Dirichlet
heat flow in a convex domain on Euclidean space,
then $F$ satisfies Assumption~\ref{ell} except for the continuity at $t=1$
by \cite{IST07}*{Theorems~1.5 and 1.6} and \cite{IST}*{Theorem~2.4 and Subsection~4.3}. 
%
 \section{Numerical experiments}
 \label{sec:numerical_experiments}
Numerical experiments demonstrate that the estimate in Theorem~\ref{thm:convergence} captures the error
\begin{equation}
\langle C,\Pi^{U}(\omega)\rangle -\inf_{\Pi\in \Pi(x,y)}\langle C,\Pi\rangle
\label{eq:error}
\end{equation}
with $\omega=(C,x,y,\varepsilon)\in \Omega$.
We compare regularization using Bregman divergence associated with the incomplete gamma function~$U_\alpha$ (subsection~\ref{powerlog}) and scaled versions of the complementary error function~$W$ (subsection~\ref{erfc}) and $U_{\mathrm{s}}$ (equation~\eqref{sym}) with regularization using the Shannon entropy~$S$ (equation~\eqref{eq:Shannon}).

All computations are performed on a computer with an Intel Core i7-8565U 1.80 GHz central processing unit (CPU), 16 GB of random-access memory (RAM), and the Microsoft Windows 11 Pro 64 bit Version~22H2 Operating System. 
All programs for implementing the method were coded and run in MATLAB R2024a for double precision floating-point arithmetic with unit roundoff $2^{-53} \simeq 1.1\!\cdot\!10^{-16}$.

The tested regularized problem~\eqref{eq:OTe} is solved by using the MATLAB built-in solver~\texttt{fmincon}.
The values of the entries of~$C \in \mathbb{R}^{3 \times 3}$ are sampled from a uniform distribution in the range~$(0,1)$ in MATLAB.
We set the tested data to $x = \begin{pmatrix} 0.1, & 0.2, & 0.7 \end{pmatrix},y = \begin{pmatrix} 0.3, & 0.4, & 0.3 \end{pmatrix} \in \mathcal{P}_3$.
Here, $\Delta_C(x,y) \simeq 4.6\!\cdot\!10^{-6}$.

\begin{figure}[ht]
    \centering
    \includegraphics[width=100mm]{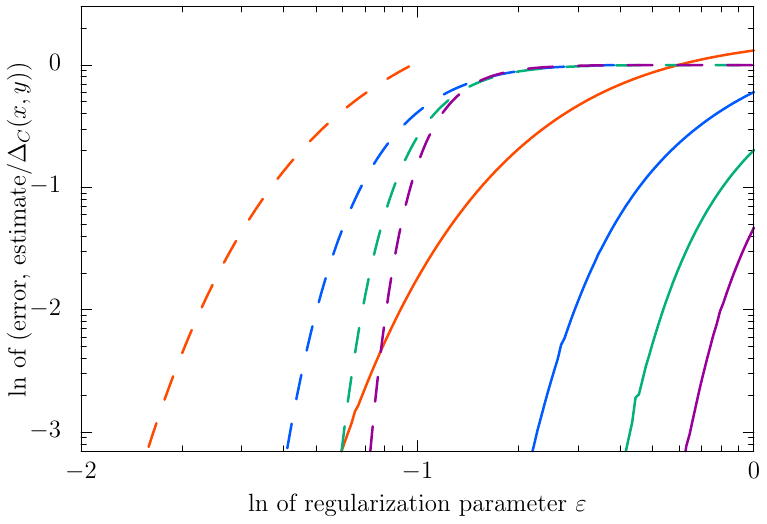}
 	\caption{Error, estimate$/\Delta_C(x,y)$ vs.\ regularization parameter~$\varepsilon$ for the incomplete gamma function~$U_\alpha$. 
  Solid line: error, dashed line: estimate, red: $\alpha=1.0$, blue: $\alpha=1/2$, green: $\alpha=1/3$, purple: $\alpha=1/4$.}
 	\label{fig:gammainc}
    \medskip

    \includegraphics[width=100mm]{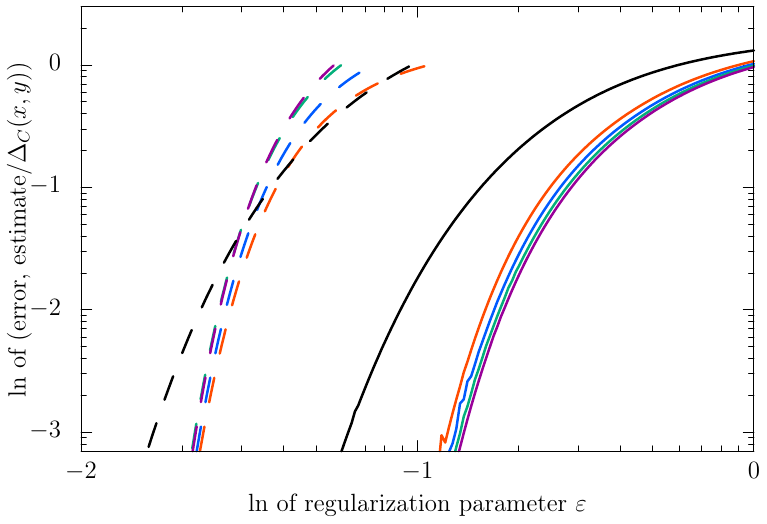}
 	\caption{Error, estimate$/\Delta_C(x,y)$ vs.\ regularization parameter~$\varepsilon$ for a scaled version of the complementary error function~$W$. 
  Solid line: error, dashed line: estimate, red: $a=2$, blue: $a=3$, green: $a=4$, purple: $a=5$, black: $U_o$.}
 	\label{fig:erfc}
 \end{figure}

\begin{figure}[t]
    \centering
    \includegraphics[width=100mm]{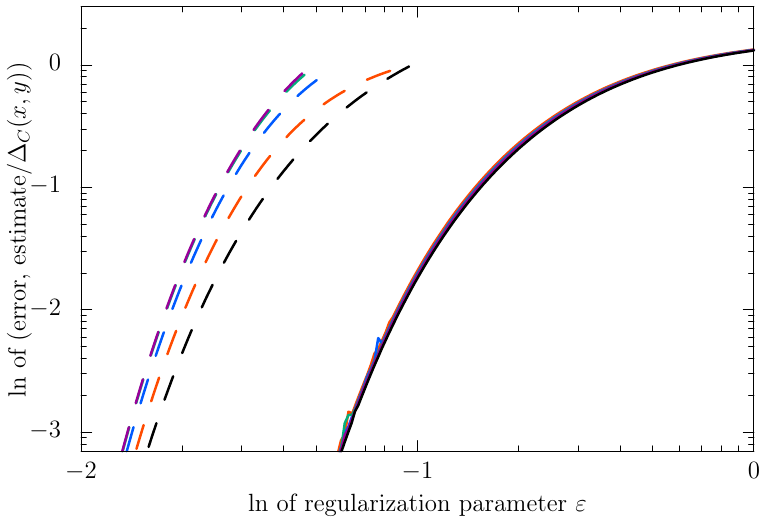}
 	\caption{Error, estimate$/\Delta_C(x,y)$ vs.\ regularization parameter~$\varepsilon$ for a scaled version of $U_{\mathrm{s}}$. 
  Solid line: error, dashed line: estimate,  red: $a=2$, blue: $a=3$, green: $a=4$, purple: $a=5$, black: $U_o$.}
 	\label{fig:U_s}
 \end{figure}

Figures~\ref{fig:gammainc}--\ref{fig:U_s} show the graphs of the natural logarithm of the 
ratio of the error~\eqref{eq:error}
and $\Delta_C(x,y)$ and its estimate
with respect to the regularization parameter~$\varepsilon$ for the incomplete gamma function~$U_\alpha$ and scaled versions of the complementary error function~$W$ and $U_{\mathrm{s}}$, respectively.
We observe that as $\varepsilon$ decreases for each value of $\alpha$ and $a$, the error decreases. 
For $U_\alpha$, as the value of $\alpha$ decreases for each value of $\varepsilon$, the error decreases.
For $U_\alpha$, as the value of $\alpha$ decreases for each value of $\varepsilon$, the estimate tends to be overestimate.
The error and its estimate for $U_\alpha$ are smaller than those for $U_o = U_{\alpha=1}$.
We next observe for $W$ that the estimate becomes smaller than that for $U_o$ for small values of $\varepsilon$.
For $W$, the error is smaller than that for $U_o$.
We finally observe for $U_\mathrm{s}$ that the estimate and error do not become smaller than those for $U_o$.
However, the decay of the estimate for $U_\mathrm{s}$ is faster than exponential.
See equation~\eqref{eu}.
As suggested from theory, with a smaller $\alpha$ for the incomplete gamma function and a larger $a$ for the complementary error function and $U_\mathrm{s}$, the error estimate becomes better.

 \section{Concluding remarks}\label{sec:condluding_remarks}

 In this paper, we considered regularization of optimal transport problems via Bregman divergence.
 We proved that the optimal value of the regularized problem converges to that of the given problem.
 More precisely, our error estimate becomes faster than exponential.
 Numerical experiments showed that the derived error estimate captures the error and regularization by a Bregman divergence outperforms that by the Kullback--Leibler divergence.

 There are several future directions subsequent to this study.
 The time complexity of our regularized problem is left open.
 It would also be interesting to extend the setting of this paper from a finite set to Euclidean space.

\begin{appendices}

 \section{Proof of Lemma~\ref{behaviorU}} 
 \label{app:proof4behaviorU}

 By the strict convexity of $U$ on $[0,1]$, $U'$ is strictly increasing on $(0,1]$ and $\lim_{h \downarrow 0}U'(h)
\in [-\infty,\infty)$ holds.
Thus, the first assertion follows.

If $U'(0)\in \mathbb{R}$,
then $\lim_{h\downarrow0}h U'(h)= 0$ holds.
%
Assume $U'(0)=-\infty$.
The Taylor expansion yields 
\[
U(r)-U(h) \geq (r-h)U'(h)
\]
for all $r,h\in(0,1]$.
By the continuity of $U$,
taking the limit as $r \downarrow0$ gives 
\begin{align*}
U(0)-U(h) \geq -h U'(h)
\end{align*}
for $h \in (0,1]$.
If $h$ is small enough, then $U'(h)<0$ by the monotonicity of $U'$ on $(0,1]$ together with $U'(0)=-\infty$.
Thus, we conclude 
\[
0=\lim_{h \downarrow0} (U(h)-U(0)) 
\leq 
 \liminf_{h\downarrow0} h U'(h)
\leq 
\limsup_{h\downarrow0} h U'(h) 
\leq 0,
\]
which leads to $\lim_{h \downarrow 0}h U'(h)= 0$.
This completes the proof of the lemma.
$\hfill\qed$

\section{Proof of Lemma~\ref{support}}
\label{app:proof4support}

For $(i,j)\in \spt(\Pi)$, we have
\[
x_i=\sum_{l=1}^{J} \pi_{il}\geq \pi_{ij}>0, \qquad
y_j=\sum_{l=1}^{I} \pi_{lj}\geq \pi_{ij}>0,
\]
which ensure that $i\in \spt(x)$ and $j\in \spt(y)$, that is, $(i,j)\in \spt(x)\times \spt(y)$.
For $(i,j)$, it turns out that 
\begin{align*}
(i,j)\in \spt(x)\times \spt(y)
\iff x_{i}>0  \text{ and }  y_j>0
\iff (i,j)\in \spt(x\otimes y).
\end{align*}
This completes the proof of the lemma.
$\hfill\qed$

\section{Proof of Lemma~\ref{lem:lower_estiamte_V}}
\label{app:proof4lower_estiamte_V}

For $r,s,t\in [0,1]$ and $r_0\in (0,1]$,
if the first inequality holds true,
then it holds that 
\begin{align}
&d_U( (1-t)r+ t s,r_0)\\
&=
U\left((1-t)r+ t s\right)
-U(r_0)
-\left[ (1-t)r+ t s- r_0 \right] U'(r_0)\\
&\geq
(1-t)U(r)+tU(s)+r U(1-t)+sU(t)-U(r_0)\\
&\quad-\left[
(1-t)r+ t s-r_0
\right]U'(r_0)\\
&=(1-t)d_U(r, r_0)
+t d_U(s, r_0)+r U(1-t)+sU(t),
\end{align}
which is the second inequality.

To show the first inequality, 
set
\[
{G}(r,s,t)
:= (1-t)U(r)+tU(s)+rU(1-t)+sU(t)- U\left( (1-t)r+ts\right)
\]
for $r,s,t \in [0,1]$.
By the continuity of $G$, 
it is enough to show 
\begin{equation}\label{enough}
\max_{r\in[0,1]}{G}(r,s,t)
\leq 0
\quad
\text{for }s,t\in (0,1).
\end{equation}
Note that 
\begin{align*}
\frac{\partial}{\partial r}{G}(r,s,t)
&=
U(1-t)+(1-t)\left(U'(r)- U'\left((1-t)r+ts\right)\right),\\
\frac{\partial^2}{\partial r^2} {G}(r,s,t)
&=(1-t)\left[ U''(r)-(1-t)U''((1-t)r+st))\right]
\end{align*}
for $r,s,t\in (0,1)$.

Fix $s,t\in (0,1)$.
Let us now show
\begin{align}
\label{eq:Phi_max}
\max_{r\in[0,1]}{G}(r,s,t)
&=
\max\left\{{G}(0,s,t),{G}(1,s,t) \right\}.
\end{align}
On one hand, for $r\in (0,1)$ with $r \leq s$, 
since $U'$ is strictly increasing on~$(0,1)$, 
we have $U'(r)-U'((1-t)r+ts)\leq 0$ and
hence $\partial_r{G} (r,s,t)\leq U(1-t)$.
Note that $U(1-t)< 0$ follows from the strict convexity of $U$ with $U(0)=U(1)=0$.
This implies that $\partial_r{G} (r,s,t) < 0 $ if $r \leq s$ and
\begin{align}\label{0s}
\max_{r\in[0,s]}{G}(r,s,t)
={G}(0,s,t), 
\quad\text{in particular, }\quad 
{G}(s,s,t)<{G}(0,s,t).
\end{align}
On the other hand, for $r\in (0,1)$ 
with $r > s$, 
we have $(1-t)r + t s < r $ and 
\[
(1-t)r U''((1-t)r + t s )
<\
 [(1-t)r + t s ] U''((1-t)r + t s )
\leq r U''(r)
\]
by $U''>0$ and the monotonicity of $r \mapsto rU''(r)$ on $(0,1)$.
This in turn yields
$\partial_r^2 {G}(r,s,t)>0$.
If $\partial_r{G} (r_0,s,t) = 0$ holds for some $r_0 \in [s,1]$,
then 
\begin{align}\label{s12}
\max_{r\in[s,1]}{G}(r,s,t)
=\max\{{G}(s,s,t), {G}(1,s,t)\}.
\end{align}
In contrast, 
if $\partial_r{G} (r,s,t) < 0$ always holds, then 
\begin{align}\label{s11}
\max_{r\in[s,1]}{G}(r,s,t)
={G}(s,s,t).
\end{align}
Summarizing the above relations~\eqref{0s}, \eqref{s12}, and \eqref{s11}, we have \eqref{eq:Phi_max}.

Since $U(0)=U(1)=0$,
a direct computation gives 
\begin{align*}
\frac{\partial^2}{\partial s^2}{G}(0,s,t)
&=
\frac{\partial^2}{\partial s^2}
\left(tU(s)+sU(t)- U\left(ts\right)\right)
=tU''(s)- t^2U''(ts)\\
&=\frac{t}{s}\left(sU''(s)- ts U''(ts) \right)\geq0
\end{align*}
for $s,t\in (0,1)$, 
where the inequality follows from the monotonicity of the function $r \mapsto r U''(r)$ on $(0,1)$.
Thus, for $t\in (0,1)$,  ${G}(0,\cdot,t)$ is convex on $[0,1]$  and
\begin{equation}\label{0st}
\max_{s\in [0,1]}{G}(0,s,t) 
=\max\{{G}(0,0,t),{G}(0,1,t) \}
=0.
\end{equation}
Next, we find 
\begin{align*}
\frac{\partial}{\partial s}{G}(1,s,t)
&=
\frac{\partial}{\partial s}
\left(tU(s)+U(1-t)+sU(t)- U(1-t+ts)\right)\\
&=
tU'(s)+U(t)- t U'\left(1-t+ts\right) 
<
U(t) 
< 0
\end{align*}
for $s,t\in (0,1)$,
where the first inequality follows from the monotonicity of $U'$ on $(0,1)$.
This leads to 
\begin{equation}\label{1st}
\max_{s\in [0,1]}{G}(1,s,t) 
={G}(1,0,t)=0
\end{equation}
for $t\in (0,1)$.
Thus, we deduce \eqref{enough} from 
\eqref{eq:Phi_max} together with~\eqref{0st} and \eqref{1st}.
This proves the lemma.
$\hfill\qed$

\section{Proof of Lemma~\ref{lem:existence_tb}}
\label{app:proof4existence_tb}

We calculate 
\[
\frac{\mathrm{d}^2}{\mathrm{d}r^2}(Dr- U(r)-U(1-r))=-U''(r)-U''(1-r)<0
\]
for $r\in (0,1)$,
consequently, 
\begin{align*}
\frac{\mathrm{d}}{\mathrm{d}r}(Dr- U(r)-U(1-r))
&>
\frac{\mathrm{d}}{\mathrm{d}r}(Dr- U(r)-U(1-r))\Bigg|_{r=R}\\
&=D-U'(R)+U'(1-R)= 0
\end{align*}
for $r\in (0,R)$.
This proves the first assertion.

We also find 
\begin{align*}
&\frac{\mathrm{d}}{\mathrm{d}r}\left( U(r)+U(1-r) -rU'(r) +r\sup_{\rho\in (0,R]}(U'(1-\rho) +\rho U''(\rho))\right)\\
&=-U'(1-r) -r U''(r)+\sup_{\rho\in (0,R]}(U'(1-\rho) +\rho U''(\rho))\\
&\geq 0
\end{align*}
for $r\in (0,R]$.
This with Lemma~\ref{behaviorU} and the assumption $U(0)=U(1)=0$ yield
\begin{align*}
&U(r)+U(1-r) -rU'(r) +r\sup_{\rho\in (0,R]}(U'(1-\rho)\\
&\geq 
\lim_{r \downarrow 0} \left(U(r)+U(1-r) -rU'(r) +r\sup_{\rho\in (0,R]}(U'(1-\rho) +\rho U''(\rho)) \right)
=0
\end{align*}
for $r\in (0,R]$.
This proves the second assertion of the lemma.
$\hfill\qed$

\section{Proof of Proposition~\ref{prop:DU}}
\label{app:proof4DU}

By a similar argument in the implication \eqref{D} to \eqref{C2}, it follows from  \eqref{K} that 
\[
\widetilde{a}^2 U''(\widetilde{a}r)=\kappa a^2 U''(ar) \quad\text{for } r \in (0,1).
\]
Let $\theta:=\widetilde{a}a^{-1}<1$.
Then, the above relation is equivalent to 
\[
\theta r U''(\theta r)=\kappa \theta^{-1} rU''(r) \quad\text{for } r \in (0,a).
\]
The monotonicity of $r \mapsto r U''(r)$ on $(0,a)$ yields $\kappa \theta^{-1} \leq 1$.
For $N\in \mathbb{N}$, it turns out that 
\begin{align*}
U'(a\theta)-U'(a\theta^{N+1})
&=\int_{a\theta^{N+1} }^{a\theta} U''(r)\mathrm{d}r
=\sum_{n=1}^{N}\int_{a\theta^{n+1}}^{a\theta^n} r U''(r) \cdot \frac1r \mathrm{d}r \\
&\leq \sum_{n=1}^{N} a\theta^n U''(a\theta^n) \cdot \frac{1}{a\theta^{n+1}} \int_{a\theta^{n+1}}^{a\theta^n} 1 \mathrm{d}r\\
&= \sum_{n=1}^{N}  \left(\kappa \theta^{-1}\right)^{n-1} a\theta U''(a\theta)\left(\theta^{-1}-1 \right)\\
&=
\widetilde{a}U''(\widetilde{a})  \left(\theta^{-1}-1 \right)
\sum_{n=1}^{N}  \left(\kappa \theta^{-1}\right)^{n-1}.
\end{align*}
If $\kappa \theta^{-1}<1$, then 
\begin{align}
\lim_{N\to \infty}(U'(a\theta)-U'(a\theta^{N+1}))
&\leq 
\widetilde{a}U''(\widetilde{a})  \left(\theta^{-1}-1 \right)
\cdot \lim_{N\to \infty} \sum_{n=1}^{N} \left(\kappa \theta^{-1}\right)^{n-1}\\
&=
\widetilde{a}U''(\widetilde{a})  \left(\theta^{-1}-1 \right)
\cdot
 \frac{1}{1-\kappa \theta^{-1}}<\infty,
\end{align}
which contradicts the condition~$\lim_{h \downarrow 0} U'(h)=-\infty$.
Hence, $\kappa \theta^{-1}=1$ and 
\[
rU''(r)=\theta rU''(\theta r) \leq rU''(r)\quad \text{for }r\in (0,a),
\]
that is, $r\mapsto r U''(r)$ is constant on $(0,a)$.
This is equivalent to that 
there exist $\mu_0, \mu_1\in \mathbb{R}$ and $\lambda>0$ such that 
$U(r)=\lambda r \log r +\mu_1 r+\mu_0$ on $(0,a]$.
This completes the proof of the proposition.
$\hfill\qed$

\end{appendices}

 \subsection*{Acknowledgements}
KM was supported in part by JSPS KAKENHI Grant Numbers~JP20K14356, \break JP21H03451, JP24K00535.
KS was supported in part by JSPS KAKENHI Grant Numbers JP22K03425, JP22K18677, JP23H00086.
AT was supported in part by JSPS KAKENHI Grant Numbers~JP19K03494, JP19H01786.
The authors are sincerely grateful to 
Maria Matveev and Shin-ichi Ohta for helpful discussion.

The authors are sincerely grateful to 
Maria Matveev and Shin-ichi Ohta for helpful discussion.
 \bibliographystyle{plain}
 \bibliography{2022_MST_OT.bib}

\begin{bibdiv}
\begin{biblist}

\bib{Ama09IEEE}{article}{
      author={Amari, Shun-Ichi},
       title={$\alpha$-divergence is unique, belonging to both $f$-divergence and {b}regman divergence classes},
        date={2009},
     journal={IEEE Trans. Inf. Theory},
      volume={55},
      number={11},
       pages={4925\ndash 4931},
}

\bib{CKLPGS22FOCS}{inproceedings}{
      author={Chen, Li},
      author={Kyng, Rasmus},
      author={Liu, Yang~P.},
      author={Peng, Richard},
      author={Gutenberg, Maximilian~Probst},
      author={Sachdeva, Sushant},
       title={Maximum flow and minimum-cost flow in almost-linear time},
        date={2022},
   booktitle={2022 {IEEE} 63rd {A}nnual {S}ymposium on {F}oundations of {C}omputer {S}cience},
   publisher={{IEEE}},
     address={Denver, CO, USA},
       pages={612\ndash 623},
}

\bib{CM94}{article}{
      author={Cominetti, R.},
      author={Mart\'{i}n, J.~San},
       title={Asymptotic analysis of the exponential penalty trajectory in linear programming},
        date={1994},
     journal={Math. Program.},
      volume={67},
       pages={169\ndash 187},
}

\bib{DMH21}{inproceedings}{
      author={Daniels, Max},
      author={Maunu, Tyler},
      author={Hand, Paul},
       title={Score-based generative neural networks for large-scale optimal transport},
        date={2021},
   booktitle={Proceedings of the {A}dvances in {N}eural {I}nformation {P}rocessing {S}ystems},
      editor={Ranzato, M.},
      editor={Beygelzimer, A.},
      editor={Dauphin, Y.},
      editor={Liang, P.S.},
      editor={Vaughan, J.~Wortman},
      volume={34},
   publisher={Curran Associates, Inc.},
     address={virtual},
       pages={12955\ndash 12965},
}

\bib{DPR18}{article}{
      author={Dessein, A.},
      author={Papadakis, N.},
      author={Rouas, J.-L.},
       title={Regularized optimal transport and the rot mover's distance},
        date={2018},
     journal={J. Mach. Learn. Res.},
      volume={19},
       pages={1\ndash 53},
}

\bib{EG}{book}{
      author={Evans, Lawrence~C.},
      author={Gariepy, Ronald~F.},
       title={Measure theory and fine properties of functions},
     edition={Revised},
      series={Textbooks in Mathematics},
   publisher={CRC Press},
     address={Boca Raton, FL},
        date={2015},
        ISBN={978-1-4822-4238-6},
}

\bib{Fan92}{article}{
      author={Fang, S.~C.},
       title={An unconstrained convex programming view of linear programming},
        date={1992},
        ISSN={0340-9422},
     journal={Z. Oper. Res.},
      volume={36},
      number={2},
       pages={149\ndash 161},
}

\bib{IST01}{article}{
      author={Ishige, Kazuhiro},
      author={Salani, Paolo},
      author={Takatsu, Aska},
       title={To logconcavity and beyond},
        date={2020},
        ISSN={0219-1997},
     journal={Commun. Contemp. Math.},
      volume={22},
      number={2},
       pages={1950009, 17},
}

\bib{IST}{article}{
      author={Ishige, Kazuhiro},
      author={Salani, Paolo},
      author={Takatsu, Asuka},
       title={Hierarchy of deformations in concavity},
        date={2022},
     journal={Info. Geo.},
}

\bib{IST07}{article}{
      author={Ishige, Kazuhiro},
      author={Salani, Paolo},
      author={Takatsu, Auka},
       title={Characterization of ${F}$-concavity preserved by the {D}irichlet heat flow},
        date={2022},
     journal={arXiv:2207.13449},
}

\bib{JST19NeurIPS}{inproceedings}{
      author={Jambulapati, Arun},
      author={Sidford, Aaron},
      author={Tian, Kevin},
       title={A direct $\tilde{O}(1/\epsilon)$ iteration parallel algorithm for optimal transport},
        date={2019},
   booktitle={Proceedings of the {A}dvances in {N}eural {I}nformation {P}rocessing {S}ystems},
      editor={Wallach, H.},
      editor={Larochelle, H.},
      editor={Beygelzimer, A.},
      editor={d\textquotesingle Alch\'{e}-Buc, F.},
      editor={Fox, E.},
      editor={Garnett, R.},
      volume={32},
   publisher={Curran Associates, Inc.},
     address={Vancouver, Canada},
}

\bib{KTM20SIMODS}{article}{
      author={Klatt, Marcel},
      author={Tameling, Carla},
      author={Munk, Axel},
       title={Empirical regularized optimal transport: {S}tatistical theory and applications},
        date={2020},
     journal={SIAM J. Math. Data Sci.},
      volume={2},
      number={2},
       pages={419\ndash 443},
}

\bib{MNPN17AAAI}{inproceedings}{
      author={Muzellec, Boris},
      author={Nock, Richard},
      author={Patrini, Giorgio},
      author={Nielsen, Frank},
       title={Tsallis regularized optimal transport and ecological inference},
        date={2017},
   booktitle={Proceedings of the {T}hirty-{F}irst {AAAI} {C}onference on {A}rtificial {I}ntelligence},
       pages={2387\ndash 2393},
}

\bib{Naudts}{book}{
      author={Naudts, Jan},
       title={Generalised {T}hermostatistics},
   publisher={Springer-Verlag},
     address={Berlin, Germany},
        date={2011},
        ISBN={978-0-85729-354-1},
         url={https://doi.org/10.1007/978-0-85729-355-8},
}

\bib{OT13}{article}{
      author={Ohta, Shin-Ichi},
      author={Takatsu, Asuka},
       title={Displacement convexity of generalized relative entropies. {II}},
        date={2013},
        ISSN={1019-8385},
     journal={Comm. Anal. Geom.},
      volume={21},
      number={4},
       pages={687\ndash 785},
}

\bib{PW09}{inproceedings}{
      author={Pele, Ofir},
      author={Werman, Michael},
       title={Fast and robust earth mover's distances},
        date={2009},
   booktitle={2009 {IEEE} 12th {I}nternational {C}onference on {C}omputer {V}ision},
       pages={460\ndash 467},
}

\bib{PeyreCuturi2019}{book}{
      author={Peyr{\'{e}}, Gabriel},
      author={Cuturi, Marco},
       title={Computational {O}ptimal {T}ransport: {W}ith {A}pplications to {D}ata {S}cience},
   publisher={Now Publishers},
     address={Delft, Netherlands},
        date={2019},
      volume={11},
      number={5-6},
}

\bib{Tak21arXiv}{article}{
      author={Takatsu, Asuka},
       title={Relaxation of optimal transport problem via strictly convex functions},
        date={2021},
     journal={arXiv:2102.07336},
}

\bib{Vil09}{book}{
      author={Villani, C\'{e}dric},
       title={Optimal {T}ransport: {O}ld and {N}ew},
   publisher={Springer-Verlag},
     address={Berlin, Germany},
        date={2009},
        ISBN={978-3-540-71049-3},
         url={https://doi.org/10.1007/978-3-540-71050-9},
}

\bib{Wee18COLT}{inproceedings}{
      author={Weed, Jonathan},
       title={An explicit analysis of the entropic penalty in linear programming},
        date={2018},
   booktitle={Proceedings of the 31st {C}onference on {L}earning {T}heory},
      editor={Bubeck, Sébastien},
      editor={Perchet, Vianney},
      editor={Rigollet, Philippe},
      series={Proceedings of Machine Learning Research},
      volume={75},
   publisher={PMLR},
     address={Stockholm, Sweden},
       pages={1841\ndash 1855},
         url={https://proceedings.mlr.press/v75/weed18a.html},
}

\end{biblist}
\end{bibdiv}
\end{document}